\documentclass[12pt,twoside,a4paper]{article}
\usepackage{amsmath,amssymb,amsthm}

\newtheorem{Thm}{Theorem}[section]
\newtheorem{Lem}[Thm]{Lemma}
\newtheorem{Pro}[Thm]{Proposition}
\newtheorem{Cor}[Thm]{Corollary}

\theoremstyle{definition}

\theoremstyle{remark}
\newtheorem{Rem}[Thm]{Remark}

\newcommand{\R}{\mathbb{R}}
\newcommand{\Z}{\mathbb{Z}}
\newcommand{\N}{\mathbb{N}}
\newcommand{\Q}{\mathbb{Q}}

\newcommand{\cO}{\mathcal{O}}

\newcommand{\cT}{\mathcal{T}}

\newcommand{\fM}{\mathfrak{M}}

\newcommand{\al}{\alpha}
\newcommand{\be}{\beta}
\newcommand{\ga}{\gamma}
\newcommand{\Ga}{\Gamma}
\newcommand{\de}{\delta}

\newcommand{\ep}{\varepsilon}

\newcommand{\si}{\sigma}

\newcommand{\ka}{\kappa}
\newcommand{\la}{\lambda}
\newcommand{\La}{\Lambda}
\renewcommand{\phi}{\varphi}

\newcommand{\es}{\emptyset}
\renewcommand{\d}{\partial}
\newcommand{\set}[2]{\{#1:\,\text{#2}\}}
\newcommand{\sm}{\setminus}
\newcommand{\sub}{\subset}

\newcommand{\ov}{\overline}
\newcommand{\wt}{\widetilde}
\newcommand{\wh}{\widehat}

\newcommand{\imm}{\operatorname{I}}
\newcommand{\hi}{\operatorname{HI}}
\newcommand{\ve}{\operatorname{VE}}
\newcommand{\vf}{\operatorname{VF}}
\newcommand{\ebd}{\operatorname{E}}
\newcommand{\fib}{\operatorname{F}}
\newcommand{\npc}{\operatorname{NPC}}
\newcommand{\sgn}{\operatorname{sign}}

\newcommand{\id}{\operatorname{id}}
\newcommand{\iso}{\operatorname{Iso}}
\newcommand{\hyp}{\operatorname{H}}

\newcommand{\hor}{\operatorname{hor}}
\newcommand{\vrt}{\operatorname{vert}}

\begin{document}

\title{Topological and geometric properties of graph-manifolds}
\author{S. Buyalo and P. Svetlov \footnote{both authors are
supported by CRDF grant RM1-2381-ST-02 and RFBR grant 02-01-00090.}}
\date{}
\maketitle

\begin{abstract} This is an exposition of results
on the existence problem of
$\pi_1$-injective
immersed and embedded surfaces in graph-manifolds, and also
of nonpositively curved metrics
on graph-manifolds, obtained by different authors. The results
are represented from a unified point of view based on the
notion of compatible cohomological classes
and some difference equation on the graph of a graph-manifold
(BKN-equation). Criteria for seven different properties of graph-manifolds
at three levels are given: at the level of compatible cohomological
classes; at the level of solutions to the BKN-equation;
in terms of spectral properties of operator invariants of a
graph-manifold.
\end{abstract}

\tableofcontents

\section{Introduction}\label{sect:introdution}

Let
$M$
be a closed three-manifold. We say that
$M$
contains a
$\pi_1$-injectively
immersed (embedded) surface, if there is an immersion
(embedding)
$g:S\to M$
of a closed surface
$S$
with nonpositive Euler characteristic which induces the
injective homomorphism
$g_\ast:\pi_1(S)\to\pi_1(M)$
of fundamental groups. We say that
$M$
contains a virtually and
$\pi_1$-injectively
embedded surface, if some finite cover of the
manifold
$M$
contains a
$\pi_1$-injectively
embedded surface. We say that a manifold
$M$
is virtually fibered over the circle, if some its finite
cover fibers over the circle with a fiber which is a closed
surface of the negative Euler characteristic. Finally, if
$M$
carries a Riemannian metric of nonpositive sectional curvature,
then we say that the manifold possesses a
$\npc$-metric.

The properties above are significantly important for the theory of
three-manifolds. One of the main conjectures of three-dimensional
topology says that any closed irreducible manifold with infinite
fundamental group contains a virtually and
$\pi_1$-injectively
embedded surface (see \cite{Sc}). According a
W.~Thurston's conjecture any closed (and even having finite volume)
hyperbolic manifold is virtually fibered over the circle
(see \cite[Conjecture~18]{T}). From geometry side, the
classification of closed manifolds admitting
$\npc$-metrics is rather nontrivial(see \cite{L}, \cite{BK2})
and by now is known up to the geometrization conjecture.

A natural problem arises, under which conditions a closed
three-manifold possesses a property from the ones listed above.
This survey is dedicated to solution of this problem in a class
$\fM$
of orientable graph-manifolds which is described below.
We restrict ourself to that class by the following reasons.
Firstly, listed above properties for the manifolds from
$\fM$
have the same nature and closely related to each other.
Secondly, as the consequence of this there is a complete solution
of our problem in the class
$\fM$.
Finally,
$\fM$
is a sufficiently general and wide class of closed
three-manifolds.

The manifolds from
$\fM$
can be described as follows. A manifold which is a trivial
$S^1$-fibration
over a surface with negative Euler characteristic and
nonempty boundary is called a block. A gluing of blocks
$M_1$
and
$M_2$
(possibly,
$M_1=M_2$)
along of boundary tori
$T_1\sub\d M_1$
and
$T_2\sub\d M_2$
is said to be regular if
$S^1$-fibers,
coming up from
$T_1$
and
$T_2$,
are not homotopic on the gluing torus. In this case,
the result of the gluing is not a block. We define
$\fM_0$
to be the class of connected, closed, orientable manifolds
glued in the regular way from blocks. Now, we define
$\fM$
as the class of orientable manifolds each of which
has a finite covering by a manifold from
$\fM_0$.
We also assume that every manifold
$M\in\fM$
contains no one-sided Klein bottle (the last condition
is a pure technical one and it is made for simplicity of
statements).

Every manifold from
$\fM$
also possesses a block structure. The blocks are Seifert fibered
spaces, i.e., foliations by circles which are in general not
fibrations (for more detail see
sect.~\ref{subsect:auxil}). It is not difficult to understand that
every manifold
$M\in\fM_0$
is irreducible, i.e., every sphere (piecewise linearly)
embedded in
$M$
bounds a ball. Thus all manifolds from
$\fM$
are also irreducible. Another invariant description of the class
$\fM$
is given below in sect.~\ref{subsubsect:gmfd}.

\subsection{Properties}

In fact, the solution of the characterization problem in
the class
$\fM$
will be described for a wider spectrum of properties than
indicated at the beginning of the article. Studying
$\pi_1$-injectively
immersed surfaces we are interested in surfaces of negative
Euler characteristic, i.e., we exclude tori (and Klein bottles)
for the reason that their immersions and embeddings are easily
to describe. It is well known that each
$\pi_1$-injectively
immersed torus in an irreducible graph-manifold, in particular,
in a manifold
$M\in\fM$,
is homotopic to a virtually embedded one and up to an isotopy
it lies in some block of
$M$
parallel to the Seifert fibration of that block
(see, for example, \cite{N3}).

Further, due to a block structure of
$M\in\fM$,
one distinguishes horizontal surfaces among surfaces immersed in
such a manifold. Namely, an immersion
$g:S\to M$
is said to be {\em horizontal,} if it is
transversal in every block to the corresponding Seifert fibration.
In that case,
$S$
has the negative Euler characteristic. It is known that every
horizontal immersion is
$\pi_1$-injective
(see \cite{RW}).

We are interested in the following properties of a manifold
$M\in\fM$:

\begin{itemize}\label{item:list}

\item[($\imm$)] $M$
contains a
$\pi_1$-injectively
immersed surface of negative Euler characteristic;

\item[($\hi$)]  $M$
contains an immersed horizontal surface;

\item[($\ebd$)] $M$
contains a
$\pi_1$-injectively
embedded surface of negative Euler characteristic;

\item[($\ve$)] $M$
contains a virtually and
$\pi_1$-injectively
embedded surface of negative Euler characteristic;

\item[($\fib$)] $M$
is fibered over the circle (with a surface of negative
Euler characteristic as a fiber);

\item[($\vf$)] $M$
is virtually fibered over the circle;

\item[($\npc$)] $M$
admits a
$\npc$-metric.

\end{itemize}

It is a remarkable and not at all obvious fact that in the class
$\fM$
all properties listed above have actually the same nature.
We explain this and give criteria for all the properties
($\imm$) -- ($\npc$)
at three levels:

\begin{itemize}

\item{} at the level of compatible cohomological classes,
see Theorem~\ref{Thm:levelcomp};

\item{} at the level of solutions to a difference equation
on the graph of a manifold (BKN-equation),
see Theorem~\ref{Thm:bknlevel};

\item{} in spectral terms of operators defined via numerical
invariants of a manifold, see sect.~\ref{sect:spec}.

\end{itemize}

The criteria in terms of compatible cohomological classes and
in terms of solutions to the BKN-equation show a common nature of
the properties and allow to find various relations between
them. On the other hand, these criteria, except some rare cases,
give no method to define using invariants of a manifold whether
it possesses or not a property from the list
($\imm$) -- ($\npc$).
To the contrary, the criteria in spectral terms are rather
constructive and allow to decide for a given manifold whether
it has or not a given property. In particular,
using these criteria we show that no one-directed arrow from
the diagram (\ref{arr:diag}) below is invertible.

The general picture of relations between the properties
($\imm$) -- ($\npc$)
can be described by the following implication diagram:

\[\begin{array}{ccccc}\label{arr:diag}
   \ebd & \Rightarrow & (\ve=\npc\cup\ebd) &\Rightarrow &\imm\\
   \Uparrow&&\Uparrow&&\Updownarrow\\
   \fib & \Rightarrow & \vf & \Rightarrow & \hi\\
        && \Uparrow && \\
        && \npc &&
   \end{array},\tag{D}\]
where an implication
${\rm A}\Rightarrow{\rm B}$
means that if a manifold from
$\fM$
possesses a property
${\rm A}$,
then it possesses a property
${\rm B}$
also. All implications in this diagram, except
$\imm\Rightarrow\hi$
and
$\npc\Rightarrow\vf$,
are almost obvious or trivial: if a manifold
$M\in\fM$
has the property
($\fib$),
then any fiber of the corresponding fibration is an embedded
horizontal surface (see \cite{WSY}); the property to be
horizontal is preserved under finite coverings.
On the other hand, the equality
$\ve=\npc\cup\ebd$
and the implications
$\imm\Rightarrow\hi$,
$\npc\Rightarrow\vf$
are not trivial (see sect.~\ref{subsubsect:npcebd},
Corollary~\ref{Cor:imm_hi} and Corollary~\ref{Cor:npc-vf} below),
as well as the fact that no one-directed arrow in the diagram
(\ref{arr:diag}) is invertible (see sect.~\ref{subsect:diag}).

The list of the properties above might be extended, for
example, by the property of a graph-manifold to be the link
of a singularity of a complex surface
(see \cite{N1}, \cite{N2}) or by the property
to have a
$\pi_1$-injectively embedded surface not homological to
sum of tori.

\subsection{Historical remarks}
Our survey is not a literature guide. Rather, we tried to give
from a unified point of view as much as possible closed
exposition of results obtained in this field by different authors.
We use many ideas of original papers referring to them in
Historical remarks scattered over the article. The survey also
contains new results as well as corrects gaps and mistakes
which we have found in publications. The new results are
the equality
($\ve$)=$(\npc)\cup(\ebd)$,
and the criteria for all the properties in terms of
compatible cohomological classes and in terms of solutions
to the BKN-equation (Theorems~\ref{Thm:levelcomp} and
\ref{Thm:bknlevel}).

The properties
($\imm$) -- ($\npc$)
were studied independently by a number of authors
applying different methods. Comparing the papers
\cite{BK2} and \cite{N3}
has shown unexpected similarity in description of the properties
($\npc$)
and
($\imm$).
Analyzing this similarity we came to the conclusion that
the properties have a common nature. In turn, this
allowed to develop a unified approach to study of them and to
discover new relations between them. Below we give a list of
papers devoted to the properties
($\imm$) -- ($\npc$)
of graph-manifolds. We do not pretend that the list is complete,
and apologize in advance to authors whose papers in this field
missed our attention.

{\bf The property
($\vf$).} In the paper \cite{LW}, it was found a topological
obstruction to this property. Namely, let
$M$
be a manifold of the class
$\fM$, $V$
be the set of its blocks (we use here notions and
notations introduced in sect.~\ref{subsubsect:gmfd}).
For each vertex
$v\in V$
of the graph
$\Ga_M$
of
$M$
one defines a number
$t_v$
as
$$t_v=|k_v|-\sum_{w\in\d v}\frac{1}{|b_w|}.$$
One of the results from \cite{LW} says if all the numbers
$t_v$, $v\in V$
are positive, then the manifold
$M$
has no finite cover fibered over the circle. In
\cite{N2}, a criterion for
($\vf$)
is proved in terms so called virtualizers. Unfortunately,
that criterion allows to decide whether a given graph-manifold
has the property
($\vf$)
or not only in some exceptional cases since it is unclear how
to find an appropriate virtualizer.
An obstacle to
($\vf$)
and a criterion similar to one from \cite{N2} is given in
\cite{WYY}. Closely related to
($\vf$)
is the paper \cite{RW}, where a criterion for a horizontal surface
to be a virtual fiber of a fibration over the circle
of some finite cover of
$M$
is given. That criterion (Lemma~\ref{Lem:up} below)
will be used in sect.
~\ref{subsect:sufficevirtprop}.
The proof of main results of \cite{Sv1},
namely the implication
$\npc\Rightarrow\vf$
and a spectral criterion for
($\vf$),
has a gap: the paper has no convincing argument showing
the existence of rational solutions to the BKN-equation.
We fill in this gap in sect. \ref{subsubsect:ratap}
and \ref{subsect:symratap}.

It is proven in \cite{WSY} that each irreducible
graph-manifold with nonempty boundary has a finite
cover which fibers over the circle. This fact reflects
a general principle saying that obstructions to the properties
we study usually disappear for graph-manifolds with
nonempty boundary.

{\bf The property
($\npc$).} For irreducible graph-manifolds all standard
obstructions to existence of
$\npc$-metrics
disappear, namely, all solvable subgroups of the fundamental
group are almost abelian and centralizers virtually split
(see, for example, \cite{CE}, \cite{GW}, \cite{LY}).
Thus quite unexpected was a first series of examples of
graph-manifolds
$M\in\fM_0$
without
$\npc$-metrics,
constructed in \cite{L}. In \cite{BK1}, it is obtained
a geometrization equation for a manifold
$M\in\fM_0$,
whose solvability is equivalent to existence of
$\npc$-metrics
on
$M$.
A criterion for
($\npc$)
in terms of a quadratic form defined via numerical invariants of
$M\in\fM_0$
is given in \cite{BK2}
(this criterion has a little flaw, see sect.~\ref{subsubsect:gap}
at the end of the survey).

It was also shown in \cite{L} that every irreducible
graph-manifold with boundary carries an
$\npc$-metric.

{\bf The properties
($\fib$)
and
($\ebd$).} The obstruction from \cite{LW} mentioned above
is also an obstruction for a manifold
$M\in\fM$
to fiber over the circle. Criteria for these properties were obtained
in \cite{N2} in terms of a reduced plumbing matrix for
$M$.
However, the criterion for
($\ebd$)
is inaccurate
(see sect.~\ref{subsect:hisrem}) being rather related to the property
of a graph-manifold to contain a
$\pi_1$-injectively embedded surface which is not
homologic to sum of tori.
A spectral criterion for
($\ebd$),
which is close to Theorem~\ref{Thm:specrebd}
below, is obtained in Thesis \cite{Sv2}.

{\bf The property
($\ve$).} Criteria for this property were obtained in \cite{N2}
in terms of virtualizers and in much more effective
spectral terms close to Theorem~\ref{Thm:specrve} below.

{\bf The properties
($\imm$)
and
($\hi$).} An equation on the graph of a manifold
$M\in\fM$,
whose solvability is equivalent to existence of a
$\pi_1$-injective
immersion
$S\to M$
of a surface of negative Euler characteristic,
was obtained in \cite{N3}. The equation turns out to be
the same as the equation from \cite{BK1}, thus we call it
the BKN-equation. In the same paper \cite{N3}, a criterion of
solvability of the BKN-equation in spectral terms was found.
The implication
$\imm\Rightarrow\hi$
is proven in \cite{Sv2}.

It was shown in \cite{RW} that every irreducible graph-manifold
with boundary contains a horizontal, properly immersed surface,
the result preceding the one from \cite{WSY} about
($\vf$)
property for graph-manifolds with boundary mentioned above.

A unified approach to study the properties
($\imm$) -- ($\npc$)
is developed in Thesis \cite{Sv2}.

\subsection{Preliminaries}\label{subsect:auxil}

Here we discuss basic notions related to Seifert fibered spaces and
graph-manifolds. For more detail on Seifert fibered spaces see \cite{Sc}.

\subsubsection{Seifert fibered spaces}\label{subsubsect:Seifert}

Consider the cylinder
$D^2\times[0,1]$
fibered by the segments
$\{y\}\times[0,1],\quad y\in D^2$.
Let
$\phi_{q,p}:D^2\to D^2$
be the rotation of the disc by
$2\pi q/p$,
where
$q$
and
$p$
are coprime integers. The factor space
$T(q,p)=D^2\times[0,1]/\{(y,0)\sim(\phi_{q,p}(y),1)\}$
with induced fiber structure by circles is called a fibered
solid torus.

Now, let
$M$
be a compact, orientable three-manifold with fiber structure by
circles (we suppose for simplicity that
$M$
contains no one-sided Klein bottle). Then
$M$
is called a Seifert fibered space, and every interior fiber
$\la\sub M\sm\d M$
has a neighborhood saturated by fibers and fiber-wise
homeomorphic to some fibered solid torus
$T(q,p)$,
where the numbers
$q=q(\la)$
and
$p=p(\la)$
depend in general on the fiber
$\la$.
The boundary
$\d M$
of a Seifert fibered space consists of tori fibered by
parallel circles. An interior fiber
$\la$
is said to be regular if
$p(\la)=1$;
otherwise it is called singular. The set
$\La$
of singular fibers is finite, all regular fibers are freely
homotopic to each other, and
$M\sm\cup_{\la\in\La}\la$
is a fibration by circles.

The base of a Seifert fibered space, i.e., the factor space
$\cO M$
of the fibers, is a 2-orbifold. Topologically
$\cO M$
is a surface which we denote by
$FM$
and call the underlying surface. The Euler characteristic
of the 2-orbifold
$\cO M$
is defined as
$$\chi(\cO M)=\chi(FM)-\sum_{\la\in\La}\left(1-\frac{1}{p(\la)}\right),$$
where
$\chi(FM)$
is the Euler characteristic of the surface
$FM$.
If
$\chi(\cO M)<0$,
then the structure of Seifert foliation on a manifold
$M$
is unique up to isotopy (see \cite{Sc}).

\medskip\noindent
{\bf Waldhausen bases.} Consider an oriented Seifert fibered
space
$M$
with non\-empty boundary and (orientable) base orbifold of
negative Euler characteristic. Let
$W$
be the set of the boundary components. We choose
an orientation of fibers in
$M$,
and let
$f_w\in H_1(T_w)$
be the homological class of an oriented fiber on the boundary
torus
$T_w\sub\d M$, $w\in W$.
Here
$H_1(T_w)=H_1(T_w;\Q)$
is the first homology group with rational coefficients
$\Q$.
In what follows, we omit the notation
$\Q$
in (co)homology groups for simplicity.

There is an induced from
$M$
orientation on the torus
$T_w$,
hence, the intersection form
$\land_w:H_1(T_w)\times H_1(T_w)\to\Q$
is well defined. Every collection of elements
$\{z_w,f_w\}_{w\in W}$,
where
$z_w\in H_1(T_{w})$,
satisfying the conditions:
$z=\oplus_{w\in W}z_w$
lies in the kernel of the inclusion homomorphism
$i_\ast:H_1(\d M)\to H_1(M)$
and
$z_w\land_wf_w=1$
for all
$w\in W$,
will be called a Waldhausen basis of the Seifert fibered space
$M$
(the elements
$z_w$, $f_w$
form a basis of
$H_1(T_w)$
for all boundary components
$w\in W$).
A Waldhausen basis always exists and is defined up to
the transformation
$z_w\mapsto z_w+n_wf_w$,
where
$n_w\in\Q$, $\sum_{w\in W}n_w=0$.
In the case when
$M$
is a trivial
$S^1$-fibration
over a surface
$F$,
the elements
$z_w$
can be taken from
$H_1(T_w;\Z)$
($f_w\in H_1(T_w;\Z)$
always), and the choice of a Waldhausen basis is
equivalent to the choice of a trivialization
$M=F\times S^1$.

\medskip\noindent
{\bf Framed Seifert fibered spaces.} Let
$M$
be a Seifert fibered space such that for every boundary
torus
$T_w\sub\d M$
an element
$c_w\in H_1(T_w;\Z)$
with
$c_w\land_wf_w\neq 0$
is fixed. Then the space
$M$
is called framed, and the collection
$C=\{c_w\}_{w\in W}$
is a framing of
$M$.
Using the following Lemma we define in sect.~\ref{subsubsect:gmfd}
charges, which are numerical
invariants of a graph-manifold playing an important role
in criteria for the properties
($\imm$) -- ($\npc$).

\begin{Lem}\label{Lem:eulnum} Let
$C=\{c_w\}$
be a framing of an oriented Seifert fibered space
$M$.
Then
$$\sum_{w\in W}\frac{1}{c_w\land_wf_w}
  \cdot(i_w)_\ast c_w=e(M,C)\cdot f,$$
where
$(i_w)_\ast:H_1(T_w)\to H_1(M)$
is the inclusion homomorphism induced by
$i_w:T_w\hookrightarrow M$,
$f\in H_1(M)$
is the homological class of a regular fiber,
$(i_w)_\ast f_w=f$
for every
$w\in W$,
and
$e(M,C)$
is a rational number.
\end{Lem}

\begin{proof} Let
$\{z_w,f_w\}_{w\in W}$
be a Waldhausen basis for
$M$.
Representing the class
$c_w$
as
$c_w=b_wz_w+d_wf_w$,
we obtain that
$b_w=c_w\land_wf_w$
and
$$\sum_{w\in W}\frac{1}{c_w\land_wf_w}\cdot(i_w)_\ast c_w=
  \left(\sum_{w\in W}\frac{d_w}{b_w}\right)f.$$
The coefficient
$e(M,C)=\sum_{w\in W}\frac{d_w}{b_w}$
is independent of the choice of the Waldhausen basis since
the left hand side is independent of it. Hence, the claim.
\end{proof}

The number
$e(M,C)$
is independent neither of the framing curves orientations nor
of the fiber orientation, and when the orientation of
$M$
is changed, it changes the sign. All of that immediately follows
from the definitions. We shall call the number
$e(M,C)$
the charge of the framed Seifert fibered space
$(M,C)$.
It differs only by sign from the relative Euler number for
$(M,C)$
introduced in \cite{LW}.

\subsubsection{Graph-manifolds}\label{subsubsect:gmfd}

At the beginning of the article, we gave a informal
description of the graph-manifold classes
$\fM_0$
and
$\fM$.
From that, it does not follow, at least immediately, that a
manifold cannot have different block decompositions. Here, we give
an invariant description of the classes
$\fM_0$
and
$\fM$,
which implies uniqueness (up to isotopy) of the (maximal) block
decomposition. Thus various characteristics related to the
decomposition are topological invariant of a graph-manifold.

Recall that the splitting of a three-manifold
$M$
along a proper surface
$\cT\sub M$
is a manifold
$M|\cT$,
which is homeomorphic to the complement
$M\sm N(\cT)$
of a regular open neighborhood
$N(\cT)$
of the surface
$\cT$.
There is a continuous projection
$\pi:M|\cT\to M$,
which is a homeomorphism
$\pi^{-1}(M\sm\cT)\to M\sm\cT$
on the complement to
$\cT$
and a 2-fold covering
$\pi^{-1}(\cT)\to\cT$
over
$\cT$,
furthermore,
$\pi^{-1}(\cT)\sub\d(M|\cT)$.

The class
$\fM$
(resp. its subclass
$\fM_0$)
consists of 3-dimensional, connected, closed,
orientable manifolds, which are not Seifert fibered spaces,
having the following property. For any manifold
$M\in\fM$
(resp.
$M\in\fM_0$)
there is a finite collection
$\cT$
of disjoint, embedded, incompressible tori in
$M$
such that
every connected component of the splitting
$M|\cT$
is a compact Seifert fibered space with orientable base
orbifold of negative Euler characteristic (resp.
a trivial
$S^1$-fibration
over a surface of negative Euler characteristic). Furthermore,
we assume that the manifolds from
$\fM$
contain no one-sided Klein bottle.

Let
$M$
be a manifold of the class
$\fM$.
It is well known (see \cite{JS}, \cite{J}) that a minimal
collection of tori
$\cT\sub M$
satisfying the definition is unique
up to isotopy. That collection is called the JSJ-surface in
$M$,
and the connected components of the splitting
$M|\cT$
are called the maximal blocks (or the vertex manifolds) of
the manifold
$M$.
Therefore, the decomposition of
$M$
into the maximal blocks and topological invariants of the
maximal blocks are topological invariants of
$M$.
We note that every irreducible graph-manifold, in particular,
a graph-manifold with
$\npc$-metric,
has a 2-fold cover from
$\fM$
and a finite cover from
$\fM_0$
(see, for example, \cite{N2}, \cite{RW}). Thus two approaches
to the definition of
$\fM$
agree.

\medskip\noindent
{\bf The graph of a graph-manifold.} To every
$M\in\fM$
one associates its graph
$\Ga=\Ga_M$
dual to the decomposition
$M=\cup_vM_v$
of the manifold into maximal blocks. In other words, the vertex set
$V$
of
$\Ga$
is the set of the maximal blocks of
$M$,
and the set
$W$
of the oriented edges of
$\Ga$
can be identified with the set of the boundary components
of all maximal blocks. Namely, an edge
$w\in W$
is directed from a vertex
$v$
to a vertex
$v'$,
if the boundary torus
$T_w\sub\d M_v$
is attached in
$M$
to the boundary torus
$T_{-w}\sub\d M_{v'}$,
where the minus sign means the reverse edge orientation. The
incompressible torus in
$M$,
which results from gluing of the boundary tori
$T_w$
and
$T_{-w}$,
we denote by
$T_{|w|}$.
The set of the edges pointing out of a vertex
$v$
is denoted by
$\d v$, $\d v\sub W$,
and if
$w\in\d v$,
then we write
$w^-=v$, $(-w)^+=v$.
The set of the nonoriented edges is denoted by
$E$.
The elements
$e\in E$
are identified with the pairs
$(w,-w)$, $w\in W$.

\medskip\noindent
{\bf Invariants of a graph-manifold.} Here we describe
numerical invariants of a manifold
$M\in\fM$
which play an important role in what follows. There are two
types of them: intersection indices and charges.

We fix an orientation of
$M$.
Then every boundary torus
$T_w$, $w\in W$,
receives the induced orientation, and the orientations of
$T_w$
and
$T_{-w}$
are opposite on
$T_{|w|}$.
Let
$a\land_w b\in\Z$
be the intersection index of integer homological classes
$a$, $b\in H_1(T_{|w|};\Z)$
w.r.t. the orientation coming up from the side of
$T_w$.
Then
$a\land_{-w}b=-a\land_w b$
since the orientations of the torus
$T_{|w|}$
coming up from different sides are opposite. The operation
$\land_w:H_1(T_{|w|};\Z)\times H_1(T_{|w|};\Z)\to\Z$
is extended by linearity to an operation
$\land_w:H_1(T_{|w|})\times H_1(T_{|w|})\to\Q$.

\medskip\noindent
{\bf The intersection index.} For every
maximal block, we fix one of two possible orientations
of its Seifert fibers. This distinguishes an element
$f_w\in H_1(T_{|w|};\Z)$
corresponding to an oriented Seifert fiber for every boundary torus
$T_w$.
The integer
$$b_w=f_{-w}\land_wf_w$$
satisfies
$b_w=b_{-w}\neq 0$
and it is called {\em the intersection index} of fibers on
the torus
$T_{|w|}$.
It changes the sign when the manifold orientation as
well as the orientation one of the fibers
$f_w$, $f_{-w}$
is changed.

\medskip\noindent
{\bf The form of intersection indices.} In general, it is not
possible to find orientations in such a way that all
intersection indices
$b_w$, $w\in W$,
would have one and the same sign. The
obstacle to this is the cohomological class
$\rho\in H^1(\Ga;\Z_2)$
(described in \cite[\S 9.6]{W}) of the cocycle
$\sgn b:W\to\Z_2$, $\sgn b(w)=\sgn b_w$
(here and in what follows we use the
multiplicative form of the group
$\Z_2$).
We call the class
$\rho=[\sgn b]$
{\em the form of intersection indices.} A 2-fold cover
always exists such that its form
$\rho$
is trivial.

In what follows, we always assume that the choice of
a Seifert fiber orientation for a maximal block
$M_v\sub M$
means also the choice of an element
$f_v\in H_1(M_v;\Z)$
corresponding to an oriented regular fiber, and of elements
$f_w\in H_1(T_{|w|};\Z)$
representing this fiber on the boundary tori
$T_{|w|}$, $w\in\d v$.
Then
$(i_w)_\ast f_w=f_v$,
where
$(i_w)_\ast:H_1(T_{|w|};\Z)\to H_1(M_v;\Z)$
is the inclusion homomorphism.

\medskip\noindent
{\bf The charge.} When fiber orientations of all maximal blocks
are fixed, each of them,
$M_v\sub M$,
has a natural framing
$C_v$.
Namely, we take as a distinguished class
$c_w\in H_1(T_{|w|};\Z)$, $w\in\d v$,
the class of oriented Seifert fibers of
the adjacent block,
$c_w=f_{-w}$.
The charge
$$k_v=e(M_v,C_v)\in\Q$$
of the framed Seifert fibered space
$(M_v,C_v)$
(see sect.~\ref{subsubsect:Seifert}) is called
{\em the charge} of the maximal block
$M_v$.
It follows from the properties of charges of framed Seifert
fibered spaces that the charges
$k_v$, $v\in V$,
are independent of the choice of Seifert fiber orientations of
maximal blocks, change the sign when the manifold orientation
is changed, and they are topological invariant of the
oriented manifold
$M$.

\medskip\noindent
{\bf The labeled graph.} One associates to every
oriented graph-manifold
$M$
of the class
$\fM$
its graph
$\Ga=\Ga(V,W)$
together with the collection of the absolute values of
the intersection indices
$|B|=\set{|b_w|\in\N}{$w\in W$}$,
the collection of the charges
$K=\set{k_v\in\Q}{$v\in V$}$
and the form of intersection indices
$\rho$.
The quadruple
$(\Ga,|B|,K,\rho)$
is called the label graph of the manifold
$M$.
A graph-manifold may not be recovered from its labeled graph.
However, the information encoded in the labeled graph is
sufficient to define whether a given graph-manifold possesses
or not any property from the list
($\imm$) -- ($\npc$).

In many situations the form of intersection indices
$\rho$
plays no role, then saying about the labeled graph we mean
the triple
$(\Ga,|B|,K)$.

\section{Compatible cohomological classes}\label{sect:compcol}

The key to a unified approach to the properties
($\imm$) -- ($\npc$)
is the notion of a compatible collection of
cohomological classes. Here we give the definition of this notion
and formulate criteria for
($\imm$) -- ($\npc$)
in its terms. These criteria
(Theorem~\ref{Thm:levelcomp}) do not look much effective if we
want to define whether a given manifold
$M\in\fM$
possesses this or that property from the list
($\imm$) -- ($\npc$).
However, the criteria at the level of the BKN-equation
(Theorem~\ref{Thm:bknlevel}) and the criteria in spectral terms
of operator invariants
(Theorems \ref{Thm:specrimmhi} -- \ref{Thm:specrnpc}),
which are sufficiently effective, are based on them.

\subsection{Motivation and definition}\label{subsect:motivdef}

To motivate the definition of a compatible collection let us
consider a
$\pi_1$-injective
immersion
$g:S\to M$
of a closed surface
$S$
of negative Euler characteristic in a manifold
$M\in\fM$
(all properties
($\imm$) -- ($\vf$)
imply the existence of such an immersion). Let
$\cT$
be the JSJ-surface in
$M$.
Then
$g$
is homotopic to an immersion such that the preimage
$g^{-1}(\cT)$
consists of a finite number of disjoint, simple, closed,
noncontractible curves on the surface
$S$,
and the connected components of the preimage of every
maximal block
$M_v$
are mapped in it either horizontally (horizontal components)
or parallel to the Seifert fibers of the block (vertical annuli),
see \cite{RW}. We assume in what follows that every
$\pi_1$-injective
immersion we consider has been put in this position.

Furthermore, one can assume that the vertical annuli of
the immersion
$g$,
lying in one and the same block, are separated into parallel pairs.
One can achieve this by taking the boundary of the collar of the
surface
$g(S)$.
We fix an orientation of
$M$
and orientations of Seifert fibers of its maximal blocks. This
induces an orientation on the horizontal components of the
intersection
$g(S)\cap M_v$
for every maximal block
$M_v\sub M$.
We orient the vertical annuli so that the parallel annuli
of every pair have opposite orientations. Therefore, a relative
class
$[g(S)\cap M_v]\in H_2(M_v,\d M_v)$
of the oriented surface
$g(S)\cap M_v$
is well defined. Let
$l_v\in H^1(M_v)$
the class dual to
$[g(S)\cap M_v]$.
The collection of classes
$\set{l_v}{$v\in V$}$
depends on the choice of Seifert fiber orientations and is
independent of the choice of orientation of
$M$.
This collection satisfies some conditions
axiomatizing which we come to the notion of compatible
cohomological classes.

For a class
$l_v\in H^1(M_v)$
and an edge
$w\in\d v$
we introduce the notation
$l_w:=i_w^\ast l_v$,
where the homomorphism
$i_w^\ast:H^1(M_v)\to H^1(T_{|w|})$
is induced by the inclusion
$i_w:T_{|w|}\hookrightarrow M_v$.

\begin{Lem}\label{Lem:motivcohom1} The collection of the
cohomological classes
$\set{l_v}{$v\in V$}$
obtained above for a
$\pi_1$-injective
immersion
$g:S\to M$
satisfies the conditions:
$|l_{-w}(f_w)|\le l_w(f_w)$
and if
$$|l_{-w}(f_w)|=l_w(f_w),\ |l_w(f_{-w})|=l_{-w}(f_{-w}),$$
then
$l_{-w}=\pm l_w$
for every edge
$w\in W$
of the graph
$\Ga_M$.
\end{Lem}

\begin{proof} Let
$c\sub g(S)\cap T_{|w|}$
be the image of a simple closed curve on the surface
$S$.
Then the curve
$c$,
being a boundary component of the surface
$g(S)\cap M_{w^-}$
as well as the surface
$g(S)\cap M_{w^+}$,
receives the orientations from both, and the
orientations not necessarily coincide. The curve
$c$
is called consistent (inconsistent) if these two orientations
agree (are opposite). Denote by
$c_w^+$
(resp.
$c_w^-$)
an element of
$H_1(T_{|w|})$,
which is the sum of homological classes of consistent (resp.
inconsistent) curves from
$g(S)\cap T_{|w|}$
oriented as the boundary of the surface
$g(S)\cap M_v$,
where
$v=w^-$
is the initial vertex of
$w$.
It follows from the definition that
$c_{-w}^+=c_w^+$
and
$c_{-w}^-=-c_w^-$.
Computing the value of the cocycle
$l_w$
on some
$x\in H_1(T_{|w|})$,
we obtain
$l_w(x)=(c_w^++c_w^-)\land_wx$.
On the other hand,
$l_{-w}(x)=(c_{-w}^++c_{-w}^-)\land_{-w}x
  =-(c_w^+-c_w^-)\land_wx$.
Thanks to the chosen orientation of the surface
$g(S)\cap M_v$,
the numbers
$a_w^+=c_w^+\land_wf_w$
and
$a_w^-=c_w^-\land_wf_w$
are nonnegative. The formulae above imply
$l_w(f_w)=a_w^++a_w^-$
and
$l_{-w}(f_w)=-(a_w^+-a_w^-)$.
Therefore, for every edge
$w\in W$
we have
$$|l_{-w}(f_w)|=|a_w^+-a_w^-|
  \le a_w^++a_w^-=l_w(f_w).$$
Suppose now that
$|l_{-w}(f_w)|=l_w(f_w)$
and
$|l_w(f_{-w})|=l_{-w}(f_{-w})$
for some edge
$w\in W$.
Then
$a_w^+\cdot a_w^-=0=a_{-w}^+\cdot a_{-w}^-$.
Assume that
$a_w^+=0$.
It means that either the set of the consistent curves on the torus
$T_{|w|}$
is empty or the consistent curves are vertical with respect to
the block
$M_v$.
In the first case, the class
$c_w^+$
vanishes, thus
$c_{-w}^+=0$
and
$a_{-w}^+=0$,
which easily implies the equality
$l_{-w}=l_w$.
In the second case, due to the choice of the orientations
on parallel pairs of the vertical annuli, there are consistent
as well as inconsistent curves on the torus
$T_{|w|}$,
which are vertical w.r.t. the block
$M_v$,
and hence horizontal w.r.t. the adjacent block
$M_{v'}$, $v'=w^+$.
Thus the numbers
$a_{-w}^+$, $a_{-w}^-$
are both nonzero. This is a contradiction with what precedes,
and the second case is excluded. Similarly, if
$a_w^-=0$,
then
$l_{-w}=-l_w$.
\end{proof}

{\em An oriented collection of cohomological classes} associates
to each choice of orientations of Seifert fibers of maximal blocks
a collection
$\set{l_v\in H^1(M_v;\R)}{$v\in V$}$
in such a way that changing the fiber orientation of some block
$M_v$
one changes the sign of the class
$l_v$
to the opposite one.
Lemma~\ref{Lem:motivcohom1} motivates the following definition.
An oriented collection of cohomological classes
$\set{l_v\in H^1(M_v;\R)}{$v\in V$}$
is said to be {\em compatible}, if not all classes are zero and
for every edge
$w\in W$
the inequality
$|l_{-w}(f_w)|\le l_w(f_w)$
holds true; furthermore, if
$|l_{-w}(f_w)|=l_w(f_w)$
and
$|l_w(f_{-w})|=l_{-w}(f_{-w})$,
then
$l_{-w}=\pm l_w$.

We stress that in this definition one talks about real
cohomological classes while Lemma~\ref{Lem:motivcohom1}
deals with rational ones.

By Lemma~\ref{Lem:motivcohom1}, the existence of compatible
cohomological classes is a necessary condition for each property
($\imm$) -- ($\vf$).
In the next section we show that this condition is also
necessary for the property
($\npc$).

\subsubsection{The case of
$\npc$-metrics}\label{subsubsect:npc}

It is well known (see, for example, the papers \cite{GW}, \cite{LY},
\cite{E}, \cite{Sch}, \cite{B1} and the book \cite{CE}), that
every
$\npc$-metric
$g$
on a graph-manifold
$M$
has a rather special form. Namely, JSJ-surface
$\cT$
for
$M$
can be chosen in such a way that every torus
$T_{|w|}$
is geodesic and flat, and the metric
$g$
locally splits along every maximal block
$M_v$,
i.e., every point
$z\in M_v$
has a neighborhood
$U_z\sub M_v$
isometric to the metric product
$F_z\times(-\ep,\ep)$,
where
$F_z$
is a surface of nonpositive Gaussian curvature.
The splitting is naturally compatible with the Seifert fiber
structure, the Seifert fibers are closed geodesics, and all
regular fibers have one and the same length
$L_v$.

We note that even if a block
$M_v$
is a trivial
$S^1$-bundle,
i.e.,
$M_v\simeq F_v\times S^1$,
the metric
$g$
may have no global splitting along
$M_v$
(see Remark~\ref{Rem:nonprod} below). However, being lifted
in the universal cover
$\wt M_v$,
it splits globally, and
$\wt M_v$
is isometric to the metric product
$A_v\times\R$,
where
$A_v$
is a surface of nonpositive Gaussian curvature with
geodesic boundary.

As usual, we fix an orientation of
$M$
and orientations of Seifert fibers of every maximal block
$M_v\sub M$.
This defines an orientation of
$\wt M_v=A_v\times\R$
and an orientation of the factor
$\R$.
The fundamental group
$\pi_1(M_v)$
acts freely on
$\wt M_v$
by isometries leaving invariant the splitting
$\wt M_v=A_v\times\R$.
The isometries leave invariant the orientation of
$M$
since the manifold
$M$
is orientable. Furthermore, they leave invariant the
orientation of the factor
$\R$
because
$M$
contains no one-sided Klein bottle. Therefore, we have
a homomorphism
$\phi_v:\pi_1(M_v)\to\R$
which assigns to every isometry
$\ga\in\pi_1(M_v)$, $\ga:A_v\times\R\to A_v\times\R$,
its shift of the factor
$\R$.
Since the group of the homomorphisms
$\pi_1(M_v)\to\R$
is canonically isomorphic to
$H^1(M_v;\R)$,
we get a class
$l_v\in H^1(M_v;\R)$
corresponding to the homomorphism
$\phi_v$.
By the choice of orientations,
$l_v(f_v)=L_v>0$
is the length of a regular fiber.

\begin{Lem}\label{Lem:motivcohom2} The described above
collection
$\set{l_v}{$v\in V$}$
of cohomological classes is compatible. Moreover, for every
edge
$w\in W$
we have
$|l_{-w}(f_w)|<l_w(f_w)$
and
$l_w(f_w)\cdot l_w(f_{-w})=
  l_{-w}(f_w)\cdot l_{-w}(f_{-w})$,
where, recall,
$l_w=i_w^\ast l_v$.
\end{Lem}

\begin{proof} The induced flat metric on the torus
$T_{|w|}$
defines a scalar product
$g_w$
on
$H_1(T_{|w|};\R)\simeq\R^2$
such that
$g_w(a,a)$
is the square of the length of a closed geodesic on
$T_{|w|}$
representing an element
$a\in H_1(T_{|w|};\Z)$.
Then the functional
$l_w:H_1(T_{|w|};\R)\to\R$
is the projection onto the line
$\R\cdot f_w$
w.r.t. this scalar product,
$g_w(a,f_w)=l_w(a)\cdot l_w(f_w)$
for every
$a\in H_1(T_{|w|};\R)$,
since
$g_w(f_w,f_w)=L_v^2=(l_w(f_w))^2$
is the square of the length of a regular fiber of
the block
$M_v$.
The scalar products
$g_w$
and
$g_{-w}$
on the torus
$T_{|w|}$
coincide, thus
$$l_{-w}(f_w)\cdot l_{-w}(f_{-w})=
  g_{-w}(f_w,f_{-w})
  =g_w(f_{-w},f_w)=l_w(f_w)\cdot l_w(f_{-w}).$$
The vectors
$f_{-w}$, $f_w\in H_1(T_{|w|};\R)$
are linearly independent, thus for their scalar product we have
$$|g_{-w}(f_w,f_{-w})|<\sqrt{g_w(f_w,f_w)}\sqrt{g_w(f_{-w},f_{-w})}=
   l_w(f_w)\cdot l_{-w}(f_{-w}),$$
which implies
$|l_{-w}(f_w)|<l_w(f_w)$.
\end{proof}

\subsection{Criteria for the properties
($\imm$) -- ($\npc$)}\label{subsect:crit}
Now, we are able to formulate criteria for the properties
($\imm$) -- ($\npc$)
in terms of compatible cohomological classes.

\begin{Thm}\label{Thm:levelcomp} Let
$M$
be a manifold of the class
$\fM$
and
$\Ga=\Ga(V,W)$
be its graph.
\begin{itemize}
\item[{\rm($\imm$)}] $M$
has the property
{\rm($\imm$)}
if and only if there is a compatible collection of
cohomological classes for
$M$;

\item[{\rm($\hi$)}] $M$
has the property
{\rm($\hi$)}
if and only if there is a compatible collection
$\set{l_v}{$v\in V$}$
of cohomological classes for
$M$
such that
$l_v(f_v)>0$
for every vertex
$v\in V$;

\item[{\rm($\ebd$)}] $M$
has the property
{\rm($\ebd$)}
if and only if there is a compatible collection
$\set{l_v}{$v\in V$}$
of cohomological classes for
$M$
such that for every edge
$w\in W$
we have: if
$l_w(f_w)\cdot l_{-w}(f_{-w})\neq 0$,
then
$l_{-w}=\pm l_w$,
where
$l_w=i_w^\ast l_v$;

\item[{\rm($\ve$)}] $M$
has the property
{\rm($\ve$)}
if and only if there is a compatible collection
$\set{l_v}{$v\in V$}$
of cohomological classes for
$M$
such that
$l_{-w}(f_w)\cdot l_{-w}(f_{-w})=
 l_w(f_w)\cdot l_w(f_{-w})$
for every edge
$w\in W$;

\item[{\rm($\fib$)}] $M$
has the property
{\rm($\fib$)}
if and only if there are a function
$\ep:V\to\{\pm 1\}$
and a compatible collection
$\set{l_v}{$v\in V$}$
of cohomological classes for
$M$
such that
$l_{-w}(f_w)=\ep_v\ep_{v'}l_w(f_w)\neq 0$
for every edge
$w\in W$,
where
$v=w^-$, $v'=w^+$;

\item[{\rm($\vf$)}] $M$
has the property
{\rm($\vf$)}
if and only if there is a compatible collection
$\set{l_v}{$v\in V$}$
of cohomological classes for
$M$
such that
$l_v(f_v)>0$
for every vertex
$v\in V$
and
$l_{-w}(f_w)\cdot l_{-w}(f_{-w})=
 l_w(f_w)\cdot l_w(f_{-w})$
for every edge
$w\in W$;

\item[{\rm($\npc$)}] $M$
has the property
{\rm($\npc$)}
if and only if there is a compatible collection
$\set{l_v}{$v\in V$}$
of cohomological classes for
$M$
such that
$|l_{-w}(f_w)|<l_w(f_w)$
and
$l_{-w}(f_w)\cdot l_{-w}(f_{-w})=
 l_w(f_w)\cdot l_w(f_{-w})$
for every edge
$w\in W$.
\end{itemize}
\end{Thm}

The proof of Theorem~\ref{Thm:levelcomp} is given in
sect.~\ref{sect:proofcr}. The implication
$\npc\Rightarrow\vf$
is an immediate corollary of Theorem~\ref{Thm:levelcomp}.

\begin{Cor}\label{Cor:npc-vf} If a closed graph-manifold
$M$
carries a
$\npc$-metric,
then it virtually fibers over the circle.
\end{Cor}

\begin{proof} As we already mentioned, the manifold
$M$
has a 2-fold cover from
$\fM$.
Thus we can assume that
$M\in\fM$.
Comparing the criteria
($\vf$)
and
($\npc$)
from Theorem~\ref{Thm:levelcomp}, we obtain that
$M$
virtually fibers over the circle.
\end{proof}

\section{BKN-equation}\label{sect:BKN}

By Theorem~\ref{Thm:levelcomp} each of the properties
($\imm$) -- ($\npc$)
is characterized by existence of corresponding compatible
cohomological classes. However, Theorem~\ref{Thm:levelcomp}
gives no method to define whether a given manifold
$M\in\fM$
possesses or not these classes. On the other hand, it turns out
that the classes considered as unknown in a sense satisfy a
difference equation on the graph of the manifold which is called
the BKN-equation. Its coefficients are charges and intersection
indices. This allows to make one step further to the solution of
the initial problem and to reformulate each of the properties
($\imm$) -- ($\npc$)
in terms of solutions to the BKN-equation.

\subsection{Deriving the BKN-equation}\label{subsect:bkn}

Let
$M$
be a manifold of the class
$\fM$, $\Ga=\Ga(V,W)$
be its graph,
$X=(\Ga,|B|,K)$
be its labeled graph. Assume that there is
a compatible collection
$\set{l_v\in H^1(M_v;\R)}{$v\in V$}$
of cohomological classes on
$M$
(we suppose that an orientation of
$M$
and orientations of Seifert fibers of the maximal blocks
are fixed). Define functions
$a:V\to\R$
and
$\ga:W\to\R$
as follows:
$a(v)=l_v(f_v)$
and
$$\ga(w)=\begin{cases}
  \sgn(b_w)l_w(f_{-w})(l_{-w}(f_{-w}))^{-1},
          &\text{if $l_{-w}(f_{-w})\neq 0$}\\
       0, &\text{if $l_{-w}(f_{-w})=0$}
       \end{cases}\ ,$$
where
$l_w=i_w^\ast l_v$, $(i_w)_\ast f_w=f_v\in H_1(M_v;\Z)$
is the class representing a regular oriented fiber of
$M_v$.
We put for brevity
$a_v=a(v)$, $\ga_w=\ga(w)$.
Note that for every edge
$w\in W$
the following equality holds
$\ga_w\cdot l_{-w}(f_{-w})=\sgn(b_w)\cdot l_w(f_{-w})$.
If
$l_{-w}(f_{-w})\neq 0$,
then the equality follows from the definition of
$\ga_w$;
if
$l_{-w}(f_{-w})=0$,
then
$l_w(f_{-w})=0$
due to compatibility of the collection
$\set{l_v}{$v\in V$}$.
Thus for every vertex
$v\in V$
we have
\begin{eqnarray*}
\sum_{w\in\d v}\frac{\ga_w}{|b_w|}a_{w^+}&=&
  \sum_{w\in\d v}\frac{\sgn(b_w)l_w(f_{-w})}{|b_w|}=
  \sum_{w\in\d v}\frac{i_w^\ast l_v(f_{-w})}{f_{-w}\land_wf_w}\\
  &=&l_v\left(\sum_{w\in\d v}
  \frac{(i_w)_\ast f_{-w}}{f_{-w}\land_wf_w}\right)=
  l_v(k_vf_v)=k_va_v.
\end{eqnarray*}
Therefore, the functions
$a$
and
$\ga$
satisfy the equation
$$k_va_v=\sum_{w\in\d v}\frac{\ga_w}{|b_w|}a_{w^+},\quad
   v\in V,$$
which is called the BKN-equation over the labeled graph
$X$.
The coefficients of the equation are numerical invariants
of the oriented manifold
$M$:
the charges
$k_v$, $v\in V$,
and the absolute values of the intersection indices
$|b_w|$, $w\in W$.
The numbers
$\ga_w$, $w\in W$,
as well as the charges are independent of the choice
of fiber orientations of the maximal blocks, and they change the
sign to the opposite one when the manifold orientation is
changed.

It also follows from the compatibility condition of
$\{l_v\}$
that the solution
$(a,\ga)$
we have obtained has the property:
$|\ga_w|\le 1$
and if
$|\ga_w|=|\ga_{-w}|=1$
for some edge
$w\in W$,
then
$\ga_{-w}=\ga_w=\pm 1$.

For each of the properties
($\imm$) -- ($\npc$)
unknown functions
$a$
and
$\ga$
have a clear geometric meaning. For instance, if
one talks about
($\npc$),
then
$a_v=L_v$
is the length of a regular fiber of the block
$M_v$,
and
$$\ga_w=\sgn(b_w)\frac{g_w(f_w,f_{-w})}
  {\sqrt{g_w(f_w,f_w)\cdot g_w(f_{-w},f_{-w})}}$$
is
$\pm$
cosine of the angle between the oriented fibers
$f_w$
and
$f_{-w}$
of adjacent fibrations on the common torus
$T_{|w|}$.
If one talks about the property
($\hi$),
then
$a_v$
is the degree of the projection of the horizontal
surface
$g(S)\cap M_v$
onto the base orbifold
$\cO M_v$
of the maximal block
$M_v$,
where
$g:S\to M$
is a horizontal immersion. In that case one holds
$|\ga_w|\le 1$
by Lemma
~\ref{Lem:motivcohom1}. A similar interpretation of the variables
$a$
and
$\ga$
takes place for every remaining property
($\imm$), ($\ebd$) -- ($\vf$).
In any case
$a_v\ge 0$
for all vertices
$v\in V$
and
$|\ga_w|\le 1$
for all edges
$w\in W$.
Admitting a conventionality, the function
$a\in\R^V$
from a solution
$(a,\ga)$
to the BKN-equation will be called the length function, and
$\ga\in\R^W$
will be the angle function.

\subsection{Criteria for the properties
($\imm$) -- ($\npc$)}\label{subsect:bkncrit}

A solution
$(a,\ga)$
to the BKN-equation is said to be compatible if
\begin{itemize}
\item{} the length function
$a$
is nonnegative and
$a\not\equiv 0$;
\item{} the angle function
$\ga$
satisfies the conditions:
$|\ga_w|\le 1$
for all
$w\in W$
and if
$|\ga_w\cdot\ga_{-w}|=1$
for some edge
$w\in W$,
then
$\ga_{-w}=\ga_w=\pm 1$;

\item{} if
$a_v=0$
for some vertex
$v\in V$,
then
$\ga_w=\ga_{-w}=0$
for all edges
$w\in\d v$.
\end{itemize}

The last condition is introduced for convenience and it
is not an essential restriction because if
$(a,\ga)$
is a solution to the BKN-equation, then for
$$\ga_w'=\begin{cases}
         \ga_w, &\text{if $a_{w^+}\cdot a_{w^-}\neq0$}\\
             0, &\text{if $a_{w^+}\cdot a_{w^-}=0$}
         \end{cases}$$
the collection
$(a,\ga')$
is a solution to the BKN-equation satisfying that condition.

A solution
$(a,\ga)$
to the BKN-equation is called symmetric if the angle
function
$\ga$
is symmetric,
$\ga_w=\ga_{-w}$
for every edge
$w\in W$.
Now we are able to formulate the criteria for the
properties
($\imm$) -- ($\npc$)
in terms of solutions to the BKN-equation.

\begin{Thm}\label{Thm:bknlevel} Let
$X=(\Ga,|B|,K,\rho)$
be the labeled graph of an oriented manifold
$M\in\fM$.
\begin{itemize}
\item[{\rm($\imm$)}] $M$
has the property
{\rm($\imm$)}
if and only if the BKN-equation over the graph
$X$
has a compatible solution;

\item[{\rm($\hi$)}] $M$
has the property
{\rm($\hi$)}
if and only if the BKN-equation has a compatible
solution
$(a,\ga)$
with positive length function
$a$,
$a_v>0$
for every vertex
$v\in V$;

\item[{\rm($\ebd$)}] $M$
has the property
{\rm($\ebd$)}
if and only if the BKN-equation has a compatible, symmetric
solution
$(a,\ga)$
such that
$\ga_w=\ga_{-w}=\pm 1$
for every edge
$w\in W$
with
$a_{w^-}\cdot a_{w^+}\neq 0$;

\item[{\rm($\ve$)}] $M$
has the property
{\rm($\ve$)}
if and only if the BKN-equation has a compatible,
symmetric solution
$(a,\ga)$;

\item[{\rm($\fib$)}] $M$
has the property
{\rm($\fib$)}
if and only if there is a compatible, symmetric solution
$(a,\ga)$
to the BKN-equation with positive length function
$a$
and the angle function
$\ga:W\to\{\pm 1\}$
whose cohomological class is the form of intersection indices,
$[\ga]=\rho\in H^1(\Ga;\Z_2)$;

\item[{\rm($\vf$)}] $M$
has the property
{\rm($\vf$)}
if and only if the BKN-equation has a compatible, symmetric
solution
$(a,\ga)$
with positive length function
$a$;

\item[{\rm($\npc$)}] $M$
has the property
{\rm($\npc$)}
if and only if the BKN-equation has a compatible, symmetric
solution
$(a,\ga)$
with positive length function
$a$
and
$\ga_w=\ga_{-w}\in(-1,1)$
for every edge
$w\in W$.
\end{itemize}
\end{Thm}

The proof of Theorem~\ref{Thm:bknlevel} is given in
sect.~\ref{sect:proofcr}.

\subsection{Historical remarks}

The BKN-equation has been obtained in \cite{BK1}, \cite{BK2}
for study the property
($\npc$),
and independently in \cite{N3} for study the properties
($\imm$), ($\hi$). A relation equivalent to the BKN-equation
for the particular case
$|\ga|\equiv 1$
has been obtained and used in \cite{WSY}. This equation
has been used in \cite{BK3} to derive a criterion for
the property
($\npc$)
with finite volume for infinite graph-manifolds. In the paper
\cite{BK4}, it was generalized to the case of an arbitrary
dimension
$n\ge 4$
for manifolds glued from maximal blocks each of which is
a trivial torus bundle over a compact surface of
negative Euler characteristic. In that case, in contrast to
the case
$n=3$,
charges are not numerical invariants, but some linear mappings.
In \cite{B2}, an interpretation of the BKN-equation
has been found as an analog of the Maxwell equations
of classical electrodynamics. The criteria
($\imm$) -- ($\vf$)
of Theorem~\ref{Thm:bknlevel} have been obtained in Thesis~\cite{Sv2}.

\section{Spectral criteria}\label{sect:spec}

Here we introduce operator invariants of a manifold
$M\in\fM$,
i.e., linear operators defined via the labeled graph
$(\Ga,|B|,K,\rho)$
(see sect.~\ref{subsubsect:gmfd}), and give criteria for
the properties
($\imm$) -- ($\npc$)
in spectral terms of these operators. This completes the
program of study the properties
($\imm$) -- ($\npc$)
for the manifolds of
$\fM$
which was indicated in Introduction.

\subsection{Operator invariants}\label{subsect:opinv}

Let
$\Ga=\Ga(V,W)$
be the graph of an oriented manifold
$M\in\fM$,
$(\Ga,|B|,K,\rho)$
be its labeled graph. If not all charges of
$M$
are zero, then we always orient
$M$
in such a way that at least one of the charges is positive.

All operator invariants we use are linear operators on
$\R^V$
symmetric w.r.t. the scalar product
$$(x,x')=\sum_{v\in V}x_vx_v',\quad x=(x_v)\in\R^V.$$
Thus we define them via corresponding quadratic forms on
$\R^V$.
Note that the space
$\R^V$
has a canonical basis which we identity with the vertex set
$V$
of the graph
$\Ga$.

We define symmetric operators
$D_M$
and
$D_M^+$
via quadratic forms
$$(D_Mx,x)=\sum_{v\in V}k_vx_v^2\quad\text{и}\quad
 (D_M^+x,x)=\sum_{v\in V}|k_v|x_v^2,$$
where
$x=(x_v)\in\R^V$.
For a symmetric function
$\la:W\to\Z_2$
(a cocycle on the graph
$\Ga$)
let us define a symmetric operator
$J_M^\la$,
$$(J_M^\la x,x)=\sum_{w\in W}\frac{\la_w}{|b_w|}x_{w^-}x_{w^+}$$
(recall that we always use the multiplicative form of the group
$\Z_2$).

\subsubsection{The operator
$A_M^+$}
We define a linear operator
$A_M^+:\R^V\to\R^V$
by
$$A_M^+=D_M^+-J_M,$$
where
$J_M=J_M^\la$
for
$\la\equiv 1$.
In terms of this operator, we formulate criteria for
the properties
($\imm$)
and
($\hi$).
As a motivation of
$A_M^+$
we prove the following Lemma.

\begin{Lem}\label{Lem:nonposeigen} If a manifold
$M\in\fM$
has one of the properties from the list
{\rm($\imm$) -- ($\npc$)},
then the operator
$A_M^+$
has a nonpositive eigenvalue.
\end{Lem}

\begin{proof} By Theorem~\ref{Thm:bknlevel}, the BKN-equation
has a compatible solution
$(a,\ga)$.
Thus
\begin{eqnarray*}
(A_M^+a,a)&=&\sum_{v\in V}\left(|k_v|a_v^2-
             a_v\sum_{w\in\d v}\frac{1}{|b_w|}a_{w^+}\right)\\
  &\le&\sum_{v\in V}a_v\sgn(k_v)\left(k_va_v-
     \sum_{w\in\d v}\frac{\ga_w}{|b_w|}a_{w^+}\right)
     =0,
\end{eqnarray*}
because
$\sgn(k_v)\ga_w\le 1$.
Since the operator
$A_M^+$
is symmetric, it has a nonpositive eigenvalue.
\end{proof}

\subsubsection{Operators
$A_\la$}
For a cocycle
$\la:W\to\Z_2$,
we put
$A_\la=D_M-J_M^\la$.
We are interested in singular properties of
$A_\la$.
Recall that a linear operator
$A:\R^V\to\R^V$
is called singular, if it has a nontrivial kernel, and
supersingular, if its kernel contains an element
$x\in\R^V$
with all coordinates
$x_v$
different from zero. We say that an operator
$A:\R^V\to\R^V$
is {\em weakly singular}, if there is a nonzero vector
$x\in\R^V$
such that
$(Ax)_v=0$
for all vertices
$v$
from the support of
$x$, $x_v\neq0$.
(The properties of supersingularity and weak singularity
depend on the choice of a basis. However, the space
$\R^V$
has the canonical basis
$V$,
and hence these properties are manifold invariants for the
corresponding operators defined below).

\begin{Lem}\label{Lem:cohomdeg} Let
$\la:W\to\Z_2$
be a symmetric function. The spectrum of the operator
$A_\la$
as well as its properties to be singular, supersingular, or
weakly singular depend only on the cohomological class
$[\la]\in H^1(\Ga;\Z_2)$.
\end{Lem}

\begin{proof} Take
$\la'=\la\cdot\si$
where
$\si$
is a coboundary, i.e.,
$\si_w=\ep_{w^-}\ep_{w^+}$
for some function
$\ep:V\to\Z_2$
and every edge
$w\in W$.
Given a vector
$x\in\R^V$,
we consider the vector
$x'\in\R^V$, $x_v'=\ep_vx_v$.
Then for every number
$\mu\in\R$
and every vertex
$v\in V$
we have
$(A_{\la'}x'-\mu x')_v=\ep_v(A_\la x-\mu x)_v$,
and the assertion follows.
\end{proof}

We formulate spectral criteria for the properties
($\ebd$)
and
($\fib$)
in terms of singularity of operators
$A_\la$.
By Lemma~\ref{Lem:cohomdeg}, one can assume that
$\la\in H^1(\Ga;\Z_2)$.

\subsubsection{The operator
$H_M$}\label{subsubsect:signcomp}

We formulate spectral criteria for the properties
($\ve$), ($\vf$)
and
($\npc$)
in terms of a symmetric operator
$H_M:\R^V\to\R^V$
which is defined below. In that cases, one has to take into
account the distribution of the charges signs on the vertex set
$V$
of the graph
$\Ga$,
what essentially complicates the definition of
$H_M$.

\medskip\noindent
{\bf The graph of sign components.} Let
$E$
be the set of nonoriented edges of the graph
$\Ga$.
Vertices
$v$, $v'\in V$
of
$\Ga$
lie in one and the same sign component if
$v=v'$
or if there is a vertex sequence
$v_1=v,\dots,v_n=v'$
such that for every
$i=1,\dots,n-1$
the vertices
$v_i$, $v_{i+1}\in V$
are connected by an edge and
$k_{v_i}k_{v_{i+1}}>0$.
This defines an equivalence relation on
$V$.
The factor set
$U=V/\sim$
is called the set of sign components, and it is the disjoint
union
$U=U_0\cup U_+\cup U_-$,
where
$U_+$($U_-$)
is the set of sign components with positive (negative)
charges, and
$U_0$
is the set of the vertices from
$V$
with zero charges. Given
$u\in U$
let
$\Ga_u$
be a connected subgraph in
$\Ga$
spanned by the vertices from the component
$u$,
i.e., the graph
$\Ga_u$
contains all edges
$e\in E$
between the vertices from
$u$.
Contracting every subgraph
$\Ga_u$, $u\in U$,
to a point we obtain a graph
$G=G(U,E_0)$,
which is called the graph of sign components of the
labeled graph
$(\Ga,|B|,K)$.
The vertex set of the graph
$G$
is
$U$,
and its set of (nonoriented) edges
$E_0$
consists of all edges
$e\in E$
connecting vertices from different sign components. If
$U\neq U_0$,
then by our assumption the set
$U_+$
is nonempty. We denote by
$p:\Ga\to G$
the canonical projection.

\medskip\noindent
{\bf Defining the function
$s:U\to\{0,\pm 1\}$.}
Recall that a graph
$G$
is said to be bipartite if its vertex set
$U$
can be represented as the disjoint union
$U=P\cup N$
in such a way that every edge of the graph connects
a vertex from
$P$
with a vertex from
$N$.
This property is equivalent to that the graph
$G$
has no cycle with odd number of edges. If the graph
$G$
is connected, then the decomposition into the parts
$P$
and
$N$
is uniquely defined up to their permutation.

We define an auxiliary function
$s:U\to\{0,\pm 1\}$,
which will be used to define the operator
$H_M$,
as follows. It the graph
$G(U,E_0)$
of the sign components of the labeled graph
$(\Ga,|B|,K)$
is not bipartite, or if
$U=U_0$,
then we put
$s(u)=0$
for all
$u\in U$.
Otherwise, we choose a decomposition
$U=P\cup N$
such that
$P\cap U_+\neq\es$.
Now, we put
$s(u)=1$,
if
$u\in P$,
and
$s(u)=-1$,
if
$u\in N$.

\medskip\noindent
{\bf Defining the operator
$H_M$.}
For every sign component
$u\in U$
we put
$W_u=W\cap p^{-1}(u)$,
and define symmetric operators
$D_u$, $J_u:\R^u\to\R^u$
via the quadratic forms
$$(D_ux,x)=s(u)\sum_{v\in u}k_vx_v^2,\quad
  (J_ux,x)=\sum_{w\in W_u}\frac{1}{|b_w|}x_{w^-}x_{w^+},$$
where
$x=(x_v)\in\R^u$
(if the set
$W_u$
is empty, then
$J_u=0$).

Now, we define a symmetric operator
$H_M:\R^V\to\R^V$
by
$$H_M=\bigoplus_{u\in U}(D_u-J_u).$$
Therefore, the edges of the graph
$\Ga$
connecting different sign components give no
contribution into the operator
$H_M$.
If all charges of
$M\in\fM$
are nonzero and have one and the same sign, then
the graph of sign components
$G$
degenerates into a point. In that case, the operators
$A_M^+$, $H_M$
coincide,
$A_M^+=H_M$.

As a motivation of
$H_M$
we prove the following Lemma.

\begin{Lem}\label{Lem:symnonposeigen} If the BKN-equation over
$M\in\fM$
has a compatible, symmetric solution
$(a,\ga)$,
then
$(H_Ma,a)\le 0$,
in particular, the operator
$H_M$
has a nonpositive eigenvalue.
\end{Lem}

\begin{proof} Using the BKN-equation, we obtain
$$(H_Ma,a)=\sum_{u\in U}\sum_{w\in W_u}
   (s(u)\ga_w-1)\frac{a_{w^-}a_{w^+}}{|b_w|}\le 0.$$
The first equality in this chain of inequalities follows
from that the sum of
$s(u)\frac{\ga_w}{|b_w|}a_{w^-}a_{w^+}$
over all edges
$w\in W$
connecting vertices from different sign components is zero,
because either
$s(u)\equiv 0$,
or each summand enters this sum twice with opposite signs
$s(u)=1,-1$,
where
$u\in U$
is the sign component from which the edge
$w$
points out. Thus the operator
$H_M$
has a nonpositive eigenvalue.
\end{proof}

\subsection{Spectral criterion for
($\imm$)
and
($\hi$)}\label{subsect:specrimmhi}
Recall that a symmetric operator
$A:\R^V\to\R^V$
is said to be positive semidefinite, if
$(Ax,x)\ge 0$
for every
$x\in\R^V$.
During the proof of Theorems~\ref{Thm:levelcomp} and
\ref{Thm:bknlevel} we obtain Corollary~\ref{Cor:imm_hi},
according to which the properties
($\imm$)
and
($\hi$)
are equivalent.

\begin{Thm}\label{Thm:specrimmhi} A manifold
$M\in\fM$
has the properties
{\rm ($\imm$)}={\rm ($\hi$)}
if and only if one of the following conditions holds:

\begin{itemize}
\item[(1)] all charges of
$M$
have one and the same sign, and the operator
$A_M^+$
is positive semidefinite and singular;
\item[(2)] the operator
$A_M^+$
has a negative eigenvalue.
\end{itemize}
\end{Thm}

This Theorem is proved in sect.~\ref{subsect:proofhi}.

\subsection{Spectral criterion for
($\ebd$)}\label{subsect:specrebd}

\begin{Thm}\label{Thm:specrebd} A manifold
$M\in\fM$
has the property
{\rm($\ebd$)}
if and only if the operator
$A_\la$
is weakly singular for some class
$\la\in H^1(\Ga;\Z_2)$.
\end{Thm}

This Theorem is proved in sect.~\ref{subsect:proofebd}.

\subsection{Spectral criterion for
($\ve$)}\label{subsect:specrve}

\begin{Thm}\label{Thm:specrve} A manifold
$M\in\fM$
has the property
{\rm($\ve$)}
if and only if the operator
$H_M$
has a nonpositive eigenvalue.
\end{Thm}

This Theorem is proved is sect.~\ref{subsect:proofve}.

\subsection{Spectral criterion for
($\fib$)}\label{subsect:specrfib}

\begin{Thm}\label{Thm:specrfib} A manifold
$M\in\fM$
has the property
{\rm($\fib$)}
if and only if the operator
$A_\rho$
is supersingular, where
$\rho\in H^1(\Ga;\Z_2)$
is the form of intersection indices.
\end{Thm}

This Theorem is proved in sect.~\ref{subsect:proofib}.

\subsection{Spectral criterion for
($\vf$)}\label{subsect:specrvf}

\begin{Thm}\label{Thm:specrvf} A manifold
$M\in\fM$
has the property
{\rm($\vf$)}
if and only if one of the following conditions holds:

\begin{itemize}
\item[(1)] the operator
$H_M$
has a negative eigenvalue;
\item[(2)] the operator
$H_M$
is positive semidefinite and supersingular.
\end{itemize}
\end{Thm}

This Theorem is proved in sect.~\ref{subsect:proofvf}.

\subsection{Spectral criterion for
($\npc$)}\label{subsect:specrnpc}

\begin{Thm}\label{Thm:specrnpc} A manifold
$M\in\fM$
has the property
{\rm($\npc$)}
if and only if the operator
$H_M$
has a negative eigenvalue, or the function
$s:U\to\{0,\pm 1\}$
from its definition is zero.
\end{Thm}

This Theorem is proved in sect.~\ref{subsect:proofnpc}.

\subsection{The implication diagram}\label{subsect:diag}

Using the spectral criteria (Theorems~\ref{Thm:specrimmhi} --
\ref{Thm:specrnpc}), we give here examples, which show that
no one-directed arrow of the diagram (\ref{arr:diag}) is invertible,
and also that there are manifolds of the class
$\fM$
possessing no property from the list
($\imm$) -- ($\npc$).

\subsubsection{Graph-manifolds
$M(\al)$}
All manifolds from the examples below belong to the class
$\fM_0$
and they all have one and the same graph
$\Ga$,
which is a linear graph with three vertices
$v_1$, $v_2$, $v_3$
and two (nonoriented) edges
$e_1$, $e_2$,
connecting the vertices
$v_1$, $v_2$
and
$v_2$, $v_3$
respectively. The maximal blocks
$M_i=M_{v_i}$, $i=1,2,3$
have the following structure. The blocks
$M_1$
and
$M_3$
are trivial bundles over the torus with hole, and the block
$M_2$
is the trivial bundle over the torus with two holes. The
gluing between the blocks
$M_1$
and
$M_2$
is one and the same for all examples, and gluings between
$M_2$
and
$M_3$
are parameterized by matrices
$\al\in GL(2,\Z)$
with determinant equals
$-1$.

To describe these gluings, we fix trivializations
$M_i=F_i\times S^1$
of the blocks, orientations of
$S^1$-factors
of all blocks and orientations of the surfaces
$F_i$, $i=1,2,3$.
For every boundary torus
$T_w$, $w\in W$,
this fixes a basis
$z_w$, $f_w$
of the group
$H_1(M_v;\Z)$,
where
$w\in\d v$.
Every element
$z_w$
represents the corresponding oriented boundary component of
$F_v$,
and
$f_w$
represents the oriented factor
$S^1$
of the block
$M_v$.
We suppose that
$e_1=(w_1,-w_1)$, $e_2=(w_2,-w_2)$,
where
$w_1\in\d v_1$, $-w_1$, $w_2\in\d v_2$
and
$-w_2\in\d v_3$.
Then the gluing
$g:T_{-w_1}\to T_{w_1}$
between the blocks
$M_1$, $M_2$
is given by
\begin{eqnarray*}
   g(z_{-w_1})&=&f_{w_1}\\
   g(f_{-w_1})&=&z_{w_1}+f_{w_1}.
\end{eqnarray*}
For
\[\al=\left[\begin{array}{cc}
      a&b\\
      c&d
\end{array}\right]\in GL(2,\Z)\]
with
$ad-bc=-1$
the gluing
$h_\al:T_{-w_2}\to T_{w_2}$
between the blocks
$M_2$
and
$M_3$
is given by
$\al$,
i.e.,
\begin{eqnarray*}
   h_\al(z_{-w_2})&=&az_{w_2}+cf_{w_2}\\
   h_\al(f_{-w_2})&=&bz_{w_2}+df_{w_2}.
\end{eqnarray*}

We denote by
$M(\al)$
the graph-manifold resulting from the gluings.

The orientations of the
$S^1$-factors
of its maximal blocks and the base surfaces define
orientations of the blocks themselves which are compatible under
the gluings because the last reverse the orientation of the
corresponding boundary tori. This induces an orientation of
the manifold
$M(\al)$.
Its intersection indices are
$$b_1:=b_{w_1}=b_{-w_1}=1,\
  b_2:=b_{w_2}=b_{-w_2}=b,$$
and the charges can be easily found using Lemma~\ref{Lem:eulnum}:
$$k_1:=k_{v_1}=1,\
  k_2:=k_{v_2}=\frac{d}{b},\
  k_3:=k_{v_3}=-\frac{a}{b}.$$

\subsubsection{A graph-manifold without
($\imm$) -- ($\npc$)}

An example of such a manifold is
$M(\al)\in\fM_0$
with the gluing
\[\al=\left[\begin{array}{cc}
             1&1\\
             4&3
            \end{array}\right].\]
The operator
$A_{\al}^+=A_{M(\al)}^+$
is given by the matrix
\[A_{\al}^+=\left[\begin{array}{rrr}
                1&-1&0\\
                -1&3&-1\\
                0&-1&1
                \end{array}\right],\]
whose eigenvalues
$1$, $2\pm\sqrt 3$,
are positive. By Lemma~\ref{Lem:nonposeigen}, the manifold
$M(\al)$
possesses no property from the list
($\imm$) -- ($\npc$).

\subsubsection{($\imm)\not\Rightarrow(\ve$)}

Here we give an example
$M(\al)$
with property
($\imm$)
and without
($\ve$).
This also shows that
($\hi)\not\Rightarrow(\vf$).
As
$\al$
we take the matrix
\[\al=\left[\begin{array}{rr}
            -1&1\\
            2&-1
            \end{array}\right].\]
The operator
$A_{\al}^+=A_{M(\al)}^+$
is given by the matrix
\[A_{\al}^+=\left[\begin{array}{rrr}
                1&-1&0\\
                -1&1&-1\\
                0&-1&1
                \end{array}\right]\]
and it has a negative eigenvalue,
$1-\sqrt 2$.
By Theorem~\ref{Thm:specrimmhi}, the manifold
$M(\al)$
has the property
($\imm$).
On the other hand, the operator
$H_{M(\al)}$
is the identity,
$H_{M(\al)}=\id$.
Thus by Theorem~\ref{Thm:specrve}, the manifold
$M(\al)$
does not have the property
($\ve$).

\subsubsection{($\ve)\not\Rightarrow(\ebd$)}

Here we give an example
$M(\al)$
with property
($\ve$)
and without
($\ebd$).
As
$\al$
we take the matrix
\[\al=\left[\begin{array}{cc}
            -3&2\\
            -1&1
            \end{array}\right].\]
The operator
$H_\al=H_{M(\al)}$
is given by the matrix
\[H_\al=\left[\begin{array}{rrr}
                1&-1&0\\
                -1&1/2&-1/2\\
                0&-1/2&3/2
                \end{array}\right]\]
and it has a negative eigenvalue because
$\det H_\al=-1$.
By Theorem
~\ref{Thm:specrve}, the manifold
$M(\al)$
has the property
($\ve$).
The graph
$\Ga$
of the manifold is simply connected, thus the group
$H^1(\Ga;\Z_2)$
is trivial, and the operator
$A_\la$
coincides with
$H_\al$
for every
$\la\in H^1(\Ga;\Z_2)$.
The principal minors of the matrix
$H_\al$
of every order are different from zero, thus the operator
$A_\la=H_\al$
is not weakly singular. By Theorem~\ref{Thm:specrebd},
the manifold
$M(\al)$
does not have the property
($\ebd$).

\subsubsection{($\ve)\not\Rightarrow(\vf$)}

Here we give an example
$M(\al)$
with property
($\ve$)
and without
($\vf$).
As
$\al$
we take the matrix
\[\al=\left[\begin{array}{cc}
            0&1\\
            1&2
            \end{array}\right].\]
The operator
$H_{\al}=H_{M(\al)}$
is given by the matrix
\[H_{\al}=\left[\begin{array}{rrr}
                1&-1&0\\
                -1&2&0\\
                0&0&0
                \end{array}\right]\]
and it is singular. By Theorem~\ref{Thm:specrve}, the manifold
$M(\al)$
has the property
($\ve$).
On the other hand, the operator
$H_{\al}$
is positive semidefinite and not supersingular because
$x_1=x_2=0$
for every vector
$x=(x_1,x_2,x_3)$
with
$H_\al x=0$.
By Theorem~\ref{Thm:specrvf}, the manifold
$M(\al)$
does not have the property
($\vf$).

\subsubsection{($\ebd)\not\Rightarrow(\fib$),
($\vf)\not\Rightarrow(\fib$),
($\vf)\not\Rightarrow(\npc$)}\label{subsubsect:example}

Here we give an example
$M(\al)$
with properties
($\ebd$), ($\vf$)
and without
($\fib$), ($\npc$).
As
$\al$
we take the matrix
\[\al=\left[\begin{array}{cc}
            0&1\\
            1&1
            \end{array}\right].\]
Since the group
$H^1(\Ga;\Z_2)$
is trivial, the operator
$A_\la$
is given by the matrix
\[A_\la=\left[\begin{array}{rrr}
                1&-1&0\\
                -1&1&-1\\
                0&-1&0
                \end{array}\right]\]
for every
$\la\in H^1(\Ga;\Z_2)$,
and it is weakly singular because
$A_\la x=(0,0,-1)$
for
$x=(1,1,0)$.
By Theorem~\ref{Thm:specrebd}, the manifold
$M(\al)$
has the property
($\ebd$).

The operator
$H_{\al}=H_{M(\al)}$
is given by the matrix
\[H_{\al}=\left[\begin{array}{rrr}
                1&-1&0\\
                -1&1&0\\
                0&0&0
                \end{array}\right]\]
and thus it is positive semidefinite and supersingular,
$H_\al x=0$
for
$x=(1,1,1)$.
By Theorem
~\ref{Thm:specrvf}, the manifold
$M(\al)$
virtually fibers over the circle. On the other hand,
the operator
$A_\rho=A_\la$
is not singular,
$\det A_\la=-1$.
By Theorem~\ref{Thm:specrfib}, the manifold
$M(\al)$
does not fibers over the circle. The operator
$H_\al$
has no negative eigenvalue, and the function
$s$
from its definition is nonzero. By Theorem~\ref{Thm:specrnpc},
the manifold
$M(\al)$
carries no
$\npc$-metric.

\subsection{Historical remarks}\label{subsect:hisrem}

A spectral criterion for the property
($\fib$)
equivalent to Theorem~\ref{Thm:specrfib}, and
a criterion for the property
($\ve$)
close to Theorem~\ref{Thm:specrve} were obtained in~\cite{N2}.
A spectral criterion for the property
($\imm$)
equivalent to Theorem~\ref{Thm:specrimmhi} (but without
the equality
($\imm$)=($\hi$))
was obtained in \cite{N3}. In the same paper, an example of
a graph-manifold from the class
$\fM$
without the property
($\imm$)
has been found, and it was shown that
($\imm)\not\Rightarrow(\ve$).
A spectral criterion for the property
($\ebd$)
was obtained in \cite{Sv2}, and for the property
($\vf$)
was formulated in \cite{Sv1}. A spectral criterion for
the property
($\npc$)
was obtained in \cite{BK2} (see, however, sect.~\ref{subsubsect:gap}).
The example~\ref{subsubsect:example} is a counterexample to
Theorem~D.1(3) from \cite{N2} since the manifold
$M(\al)$
from that example has the property
($\ebd$),
and at the same time the matrix
$A_\la$
(``decomposition matrix'' in terms of \cite{N2}) is not
singular.

\section{Proof of Theorems~\ref{Thm:levelcomp} and
\ref{Thm:bknlevel}}\label{sect:proofcr}

We prove Theorems~\ref{Thm:levelcomp} and \ref{Thm:bknlevel}
simultaneously. The first step of the proof is the following
simple but important Lemma about extension.

\begin{Lem}\label{Lem:extension} Assume that for a block
$M_v$
of an oriented manifold
$M\in\fM$
an orientation of Seifert fibers is fixed and a Waldhausen basis
$\{z_w,f_w\}_{w\in\d v}$
is chosen. Given a collection of cohomological classes
$\set{l_w\in H^1(T_{|w|};\R)}{$w\in\d v$}$,
there is a class
$l_v\in H^1(M_v;\R)$
such that
$i_w^\ast l_v=l_w$
for all
$w\in\d v$
if and only if the following conditions hold
\begin{itemize}
\item[(1)] $l_w(f_w)=a_v$
is independent of
$w\in\d v$;
\item[(2)] $\sum_{w\in\d v}l_w(z_w)=0$.
\end{itemize}
\end{Lem}

\begin{proof} The conditions (1) and (2) are necessary since
for every class
$l_v\in H^1(M_v;\R)$
with
$i_w^\ast l_v=l_w$, $w\in\d v$,
we have: the number
$l_w(f_w)=l_v(f_v)$
is independent of
$w\in\d v$
and
$\sum_{w\in\d v}l_w(z_w)=l_v(\oplus_{w\in\d v}(i_w)_\ast z_w)=0$.

Conversely, assume that a collection of cohomological classes
$\{l_w\}$
satisfies the conditions (1) and (2). We show that the class
$m_v=\bigoplus_{w\in\d v}l_w\in H^1(\d M_v;\R)$
lies in the image of the homomorphism
$i^\ast:H^1(M_v;\R)\to H^1(\d M_v;\R)$.
To this end, it suffices to show, due to the exact cohomological
sequence of the pair
$(M_v,\d M_v)$,
that
$m_v(\d S)=0$
for every relative class
$S\in H_2(M_v,\d M_v;\R)$,
where
$\d:H_2(M_v,\d M_v;\R)\to H_1(\d M_v;\R)$
is the boundary homomorphism. Representing the class
$\d S$
as
$\d S=\oplus_{w\in\d v}c_w$
and decomposing the element
$c_w$
over the basis
$z_w$, $f_w$,
we obtain
$c_w=\al z_w+\be_wf_w$,
where the coefficient
$\al$
is independent of
$w\in\d v$,
because the classes
$\d S$
and
$z=\oplus_{w\in\d v}z_w$
from
$H_1(\d M_v;\R)$
lie in the kernel of the homomorphism
$i_\ast$.
Thus
$$m_v(\d S)=\sum_{w\in\d v}l_w(c_w)=(\sum_{w\in\d v}\be_w)\cdot a_v.$$
On the other hand,
$$0=i_\ast(\d S)=\sum_{w\in\d v}(i_w)_\ast(c_w)=
    (\sum_{w\in\d v}\be_w)\cdot f_v.$$
Consequently,
$m_v(\d S)=0$
and the class
$m_v$
lies in the image of the homomorphism
$i^\ast$.
Then any class
$l_v\in H^1(M_v;\R)$
with
$i^\ast l_v=m_v$
is a required one.
\end{proof}

\subsection{Theorems~\ref{Thm:levelcomp}
and \ref{Thm:bknlevel} are equivalent}

As a corollary of Lemma~\ref{Lem:extension} we obtain that
Theorems~\ref{Thm:levelcomp} and \ref{Thm:bknlevel}
are equivalent in the following sense.

\begin{Pro}\label{Pro:compequivbkn} Let
$M\in\fM$
an oriented manifold. For each property
{\rm(A)} from the list
{\rm($\imm$) -- ($\npc$)}
we have: a compatible collection of cohomological classes on
$M$
satisfying the condition {\rm(A)} of Theorem~\ref{Thm:levelcomp}
exists if and only if the BKN-equation over
$M$
has a compatible solution satisfying the condition {\rm(A)}
of Theorem~\ref{Thm:bknlevel}.
\end{Pro}

\begin{proof} Let
$\set{l_v\in H^1(M_v;\R)}{$v\in V$}$
be compatible cohomological classes satisfying the condition~(A)
of Theorem~\ref{Thm:levelcomp} (we assume that fiber orientations
of the maximal blocks are fixed). Then the functions
$a:V\to\R$, $a(v)=l_v(f_v)$,
and
$\ga:W\to\R$,
$$\ga(w)=\begin{cases}
  \sgn(b_w)l_w(f_{-w})(l_{-w}(f_{-w}))^{-1},
       &\text{if $l_w(f_w)\cdot l_{-w}(f_{-w})\neq 0$}\\
       0, &\text{if $l_w(f_w)\cdot l_{-w}(f_{-w})=0$}
       \end{cases}\ ,$$
where
$l_w=i_w^\ast l_v$,
form a compatible solution to the BKN-equation
(see sect.~\ref{subsect:bkn}). It is straightforward
to check that the solution satisfies the condition (A)
of Theorem~\ref{Thm:bknlevel}.

Conversely, let
$(a,\ga)$
be a compatible solution to the BKN-equation satisfying
the condition (A) of Theorem~\ref{Thm:bknlevel}.
For every vertex
$v\in V$,
choose a Waldhausen basis
$\{z_w,f_w\}_{w\in\d v}$
of
$M_v$.
Then for the coefficients of the decomposition
$f_{-w}=b_wz_w+d_wf_w$
we have:
$b_w$
is the intersection index for the edge
$w\in\d v$,
and
$k_v=\sum_{w\in\d v}\frac{d_w}{b_w}$
is the charge of the vertex
$v$
(see Lemma~\ref{Lem:eulnum} and sect. \ref{subsubsect:gmfd}).
For an edge
$w\in\d v$,
define a class
$l_w\in H^1(T_{|w|};\R)$
by
$l_w(f_w):=a_v$
and
$l_w(z_w):=\frac{\ga_w}{|b_w|}a_{w^+}-\frac{d_w}{b_w}a_v$.
Then
$$\sum_{w\in\d v}l_w(z_w)=
  \sum_{w\in\d v}\frac{\ga_w}{|b_w|}a_{w^+}-k_va_v=0,$$
thus the collection of the classes
$\set{l_w}{$w\in\d v$}$
satisfies the conditions of Lemma~\ref{Lem:extension}.
By that Lemma, there is a class
$l_v\in H^1(M_v;\R)$
for which
$i_w^\ast l_v=l_w$
for all
$w\in\d v$.
Since
$l_w(f_{-w})=\sgn(b_w)\ga_wa_{w^+}$
and
$|\ga_w|\le 1$,
for the vertex
$v'=w^+$
we have
$|l_w(f_{-w})|\le a_{v'}=l_{-w}(f_{-w})$.
Assume that the equalities
$|l_{-w}(f_w)|=l_w(f_w)$
and
$|l_w(f_{-w})|=l_{-w}(f_{-w})$
hold. If
$l_w(f_w)=0=l_{-w}(f_{-w})$,
then
$l_w=l_{-w}=0$.
Otherwise
$l_w(f_w)\cdot l_{-w}(f_{-w})\neq 0$,
because the solution is compatible. Then
$|\ga_w|=|\ga_{-w}|=1$,
and hence
$\ga_{-w}=\ga_w$,
because, we recall,
$\ga_w\cdot\ga_{-w}\neq -1$.
It easily follows that
$l_{-w}=\pm l_w$.
Therefore, the cohomological classes
$\set{l_v}{$v\in V$}$
are compatible. It is straightforward to check that these
classes satisfy the condition (A) of Theorem
~\ref{Thm:levelcomp}.
\end{proof}
Proposition~\ref{Pro:compequivbkn} gives an additional flexibility
for the proof of Theorems~\ref{Thm:levelcomp} and \ref{Thm:bknlevel}:
it is convenient to use compatible cohomological classes
in some cases, and compatible solutions to BKN-equations
in other ones.

\subsection{Local conditions of extension}\label{subsect:locext}

Here we establish necessary and
sufficient boundary conditions for every maximal block
to have an extension to a
$\npc$-metric
on it (Lemma~\ref{Lem:locextnpc}), and to have an extension to
a horizontal immersion in it
(Lemma~\ref{Lem:locextimm}). These conditions are
an essential part of the proof of Theorems~\ref{Thm:levelcomp}
and \ref{Thm:bknlevel}. Though Lemmas~\ref{Lem:locextnpc}
and \ref{Lem:locextimm}
are formulated differently and their proofs are technically
different, there is an explicit analogy between them.

\begin{Lem}\label{Lem:locextnpc} Assume that for a maximal block
$M_v$
of a manifold
$M\in\fM$
an orientation of its Seifert fibers is fixed and a flat metric
on the boundary
$\d M_v$
is given, i.e., for every
$w\in\d v$
a positive definite quadratic form
$g_w$
on
$H_1(T_{|w|};\R)$
is given, and
\begin{itemize}
\item[(1)] $g_w(f_w,f_w)=a_v^2>0$
is independent of
$w\in\d v$;
\item[(2)] $\sum_{w\in\d v}g(z_w,f_w)=0$,
\end{itemize}
where
$\{z_w,f_w\}_{w\in\d v}$
is a Waldhausen basis of the block
$M_v$.
Then there is a
$\npc$-metric
$g_v$
on
$M_v$
which extends the given one on the boundary
$\d M_v$,
i.e., every boundary torus
$T_{|w|}$, $w\in\d v$,
is flat and geodesic w.r.t.
$g_v$,
and the metric
$g_v$
induces on
$T_{|w|}$
the metric
$g_w$.
\end{Lem}

\begin{Rem} Condition~(2) is independent of the choice
of a Waldhausen basis; Conditions (1) and (2) are also
necessary for the existence of a
$\npc$-metric
on
$M_v$
with flat geodesic boundary (see Lemmas~\ref{Lem:motivcohom2}
and \ref{Lem:extension}).
\end{Rem}

\begin{proof}[Sketch of the proof of Lemma~\ref{Lem:locextnpc}]
Given
$w\in\d v$,
define
$l_w\in H^1(T_{|w|};\R)$
by
$l_w(a)=\frac{g_w(a,f_w)}{\sqrt{g_w(f_w,f_w)}}$,
$a\in H_1(T_{|w|};\R)$.
In other words, the class
$l_w$
acts as the projection onto
$\R\cdot f_w$
w.r.t. the scalar product
$g_w$
(cf. the proof of Lemma~\ref{Lem:motivcohom2}).
It follows from (1) and (2) that the collection
$\set{l_w}{$w\in\d v$}$
satisfies the conditions of Lemma~\ref{Lem:extension}
about extension, and thus there is a class
$l_v\in H^1(M_v;\R)$,
which induces a collection
$\{l_w\}$, $i_w^\ast l_v=l_w$,
on the boundary
$\d M_v$.
The class
$l_v$
defines a homomorphism
$\phi_v:\pi_1(M_v)\to\R$
using which we shall construct a representation
$\psi_v:\pi_1(M_v)\to\iso(\hyp^2\times\R)$
of the fundamental group of the block
$M_v$
in the isometry group of the space
$\hyp^2\times\R$.

Now, we assume for simplicity that the block
$M_v$
is the trivial
$S^1$-bundle
over a compact surface
$F_v$
of negative Euler characteristic. Fix a trivialization
$M_v=F_v\times S^1$
and an orientation of the surface
$F_v$,
what is equivalent to the choice of a Waldhausen basis
$\{z_w,f_w\}_{w\in\d v}$.
In this case, the elements
$z_w$, $w\in\d v$
represent corresponding oriented components of the boundary
$\d F_v$.

Recall that, since the Euler characteristic of
$F_v$
is negative, given any positive numbers
$L_w$, $w\in\d v$,
there is a metric of constant curvature
$-1$
with geodesic boundary on the surface
$F_v$
such that the length of
$w$-component
equals
$L_w$.

Now, as
$L_w$, $w\in\d v$,
we take the length of
$z_w$
projected on the orthogonal to
$f_w$
direction w.r.t. the metric
$g_w$.
This defines as above a hyperbolic metric on the surface
$F_v$
and hence a representation
$\eta_v:\pi_1(F_v)\to\iso(A_v)$
of the group
$\pi_1(F_v)$
in the isometry group of the universal cover
$A_v\sub\hyp^2$
of the surface
$F_v$.

We define the representation
$\psi_v:\pi_1(M_v)\to\iso(A_v\times\R)$
as
$\psi_v(\ga)=(\eta_v\circ\pi_v(\ga),\phi_v(\ga))$
for
$\ga\in\pi_1(M_v)$,
where
$\pi_v:\pi_1(M_v)\to\pi_1(F_v)$
is the projection homomorphism onto the first factor.
It is easily to check that its image acts discretely and freely
on the metric product
$A_v\times\R$,
and this gives a required
$\npc$-metric
on the block
$M_v$.

In the general case the argument is similar: one needs to use
that there is a hyperbolic structure on the orbifold
$\cO_v=\cO M_v$
with prescribed lengths
$L_w$
of the boundary components, what gives a representation
$\eta_v:\pi_1(\cO_v)\to\iso(A_v)$,
and then to take the homomorphism
$\pi_v$
from the exact sequence
$$1\longrightarrow\Z\cdot f_v\longrightarrow\pi_1(M_v)
   \stackrel{\pi_v}{\longrightarrow}\pi_1(\cO_v)\longrightarrow 1$$
of the Seifert bundle for
$M_v$.
\end{proof}

\begin{Rem}\label{Rem:nonprod} The
$\npc$-metric
constructed above on the trivial bundle
$M_v\simeq F_v\times S^1$
typically is not a metric product: the holonomy of the circle
$S^1$
along some noncontractible loops in
$F_v$
can be nontrivial. A metric is the product if it is possible to find
a trivialization
$z\in H_1(\d M_v;\Z)$, $z=\oplus_{w\in\d v}z_w$,
such that the element
$z_w$
is orthogonal to
$f_w$
w.r.t. the metric
$g_w$
for all
$w\in\d v$.
\end{Rem}

\begin{Lem}\label{Lem:locextimm} Assume that for an oriented
block
$M_v$
of a manifold
$M\in\fM$
an orientation of its Seifert fibers is fixed and for every
$w\in\d v$
two elements
$c_w^+$, $c_w^-\in H_1(T_{|w|};\Z)$
are given such that
$\al_w^\pm=c_w^\pm\land_wf_w\ge 0$,
the class
$c_w=c_w^++c_w^-$
is even, i.e.,
$\frac{1}{2}c_w\in H_1(T_{|w|};\Z)$,
and the following conditions hold
\begin{itemize}
\item[(1)]
$c_w\land_wf_w=a_v>0$
is independent of
$w\in\d v$;
\item[(2)] $\sum_{w\in\d v}z_w\land_wc_w=0$,
\end{itemize}
where
$\{z_w,f_w\}_{w\in\d v}$
is a Waldhausen basis of the block
$M_v$.
Then there is an integer
$d(M_v)\ge 1$
such that for any integer multiple
$d\ge 1$
of
$d(M_v)$
there is a horizontal immersion
$g_v:S_v\to M_v$
of a compact surface
$S_v$
with boundary
$\d S_v=\cup_{w\in\d v}(\ga_w^+\cup\ga_w^-)$,
the union of the connected components, such that
$[g_v(\ga_w^\pm)]=d\cdot c_w^\pm$
for all
$w\in\d v$,
for which
$\al_w^+\cdot\al_w^-\neq 0$.
In the case
$\al_w^+\cdot\al_w^-=0$
for some
$w\in\d v$,
one holds
$[g_v(\ga_w^+)]=[g_v(\ga_w^-)]=\frac{d}{2}\cdot c_w$.
\end{Lem}

\begin{Rem} Condition (2) is independent of the choice
of a Waldhausen basis; Conditions (1) and (2) are also
necessary for the existence of an immersion
$g_v:S_v\to M_v$
with prescribed behavior on the boundary (see
Lemmas~\ref{Lem:motivcohom1}
and \ref{Lem:extension}).
\end{Rem}

The proof of Lemma~\ref{Lem:locextnpc} was based on the existence
of a hyperbolic metric with prescribed lengths of the boundary
components on a given surface (orbifold) of negative
Euler characteristic. The proof of Lemma~\ref{Lem:locextimm}
is based on a similar existence property for coverings of
a given surface of negative Euler characteristic with
prescribed behavior on the boundary. Here is the
precise statement.

\begin{Lem}\label{Lem:presribdeg} Let
$F$
be a compact orientable surface of positive genus with
nonempty boundary,
$W$
be the set of its boundary components. Assume that every
boundary component
$w\in W$
is covered by a collection
$S_w$
consisting of
$n_w\ge 1$
circles with multiplicities
$a_w^s\ge 1$, $s\in S_w$,
where the sum
$\sum_{s\in S_w}a_w^s=a$
is independent of
$w$
and the number of the circles
$\sum_wn_w$
has the same parity as
$a\cdot|W|$.
Then the covering of the boundary
$\d F$
can be extended to a covering
$S\to F$
of multiplicity
$a$
by a connected surface
$S$.
\end{Lem}

\begin{proof} It suffices to find a homomorphism
of the group
$\pi_1(F)$
into the permutation group of a fiber (having cardinality
$a$)
with transitive image, which maps the boundary components
into permutations with the prescribed cycle structure.
Due to the parity condition, the product of the boundary
cycles is an even permutation, and thus by a result from
\cite{JLLP}, it is the commutator of an
$n$-cycle
and an involution. Now, since the genus of
$F$
is positive, the required homomorphism does exist.
\end{proof}

\begin{proof}[Proof of Lemma~\ref{Lem:locextimm}] First, we assume
for simplicity that
$M_v$
is the trivial
$S^1$-bundle
and that the trivialization
$M_v=F_v\times S^1$
corresponds to the chosen Waldhausen basis. Taking if
necessary a covering of
$F_v$
with multiplicity
$d(M_v):=|\d v|$
and with the same number of the boundary components, we can
assume that the genus of
$F_v$
is positive.

In the decompositions
$c_w^\pm=\al_w^\pm z_w+\be_w^\pm f_w$,
the integer
$\al_w^\pm=c_w^\pm\land_w f_w$
are nonnegative,
$\al_w^++\al_w^-=a_v>0$
is even for all
$w\in\d v$,
and
$\sum_{w\in\d v}(\be_w^++\be_w^-)=0$.
We construct a required immersion
$g_v:S_v\to F_v\times S^1$
as
$g_v=(g_v^{\hor},g_v^{\vrt})$.
Here
$g_v^{\hor}:S_v\to F_v$
is a covering of the degree
$a_v$,
which on a boundary component
$\ga_w^\pm\sub\d S_v$
covers the corresponding component of
$\d F_v$
with degree
$\al_w^{\pm}$,
if both
$\al_w^+$, $\al_w^-$
are nonzero, and both components
$\ga_w^\pm$
cover the corresponding component of
$\d F_v$
with degree
$a_v/2$,
if
$\al_w^+\cdot\al_w^-=0$.
Lemma~\ref{Lem:presribdeg} provides such a covering
with a connected surface
$S_v$.
The parity condition is fulfilled since the number of
the boundary components of
$S_v$
is twice of those of
$F_v$,
and
$a_v$
is even.

We require that the map
$g_v^{\vrt}:S_v\to S^1$
has the degree
$\be_w^\pm$
on the boundary component
$\ga_w^\pm$,
if
$\al_w^+\cdot\al_w^-\neq0$,
and the degree
$\frac{1}{2}(\be_w^++\be_w^-)$
on both boundary components
$\ga_w^+$, $\ga_w^-$,
if
$\al_w^+\cdot\al_w^-=0$
($\be_w^++\be_w^-$
is even because the class
$c_w$
is even).
Due to the connectedness of
$S_v$
and the condition
$\sum_{w\in\d v}(\be_w^++\be_w^-)=0$,
such a map does exist since
$[S_v,S^1]=H^1(S_v;\Z)$.
Now, the immersion
$g_v$
has the required properties by the construction.

In the general case, we remove from
$M_v$
the interiors of fibered solid tori (with disjoint closures),
which are boundaries of the singular fibers, and denote by
$\d^sv$
the set of the new boundary components of the complement
$M_v'$
to the removed pieces. The manifold
$M_v'$
is a trivial
$S^1$-bundle
over a compact orientable surface
$F_v'$
of negative Euler characteristic, furthermore, for every
torus
$T_w\sub\d M_v'$, $w\in\d^sv$,
the class
$c_w\in H_1(T_w;\Z)$
of its meridian, i.e., a simple, closed curve bounding a disc
in the corresponding removed solid tori,
is fixed. We assume that the orientations of the meridians are
chosen such that
$\al_w:=c_w\land_wf_w>0$
for all
$w\in\d^sv$.
We put
$a_v'=a_v\cdot\prod_{w\in\d^sv}\al_w$,
$\d^\ast v=\d v\cup\d^sv$,
and for
$w\in\d^\ast v$
let
$d_w=a_v'/\al_w\in\Z$,
where
$\al_w=a_v$
for
$w\in\d v$.
Then all integer classes
$c_w'=d_wc_w$, $w\in\d^\ast v$
have one and the same intersection index with
regular fiber
$f_w$,
$c_w'\land_wf_w=a_v'$.
Fix a trivialization
$M_v'=F_v'\times S^1$.
Then for the corresponding Waldhausen basis
$\{z_w,f_w\}_{w\in\d^\ast v}$
we have
$\sum_{w\in\d^\ast v}z_w\land_wc_w'=0$.

For the collection of the classes
$c_w'^{\pm}=d\cdot c_w^\pm$, $w\in\d v$,
where
$d=d(M_v):=a_v'/a_v$
(is multiplied by
$|\d^\ast v|$
if necessarily for positivity of the genus), and
$c_w'$, $w\in\d^sv$,
we construct an immersion
$g_v':S_v'\to M_v'$,
$g_v'=(g_v^{\hor},g_v^{\vrt})$
as above with only distinction that every
boundary component
$w\in\d^sv$
of the surface
$F_v'$
is covered by
$d_w$
boundary components of the surface
$S_v'$
with degree
$\al_w$
on each one under the covering
$g_v^{\hor}$.
The parity condition is reduced to that the number
$\sum_{w\in\d^sv}d_w$
has to be even, and it is fulfilled since this sum is
an integer multiple of even
$a_v$.
For the map
$g_v^{\vrt}:S_v'\to S^1$,
an additional condition is that it has the degree
$z_w\land_wc_w'$
on each of
$d_w$
corresponding boundary components of
$S_v'$
for every
$w\in\d^sv$.

One can assume that the image
$g_v'(S_v')\sub M_v'$
meets every boundary component
$w\in\d^sv$
over a collection of
$d_w$
curves parallel to the meridian. Attaching the discs to these
curves in the corresponding solid tori, we obtain
the required immersion
$g_v$
for
$d=d(M_v)$.
Now, the case of an arbitrary multiple
$d$
of
$d(M_v)$
easily follows.
\end{proof}

\subsubsection{Historical remarks}

The statement close to Lemma~\ref{Lem:locextnpc}
is given in \cite{L}. Lemma~\ref{Lem:locextimm} is a
particular case of Lemma~3.1 from \cite{N3}.
Lemma~\ref{Lem:presribdeg} and its proof are taken from \cite{N3}.

\subsection{Properties ($\imm$), ($\hi$), ($\ebd$), ($\fib$), ($\npc$)}

\subsubsection{Existence of compatible cohomological classes}

Here we show that any of the properties
($\imm$), ($\hi$), ($\ebd$), ($\fib$), ($\npc$)
for a manifold
$M\in\fM$
implies the existence of a compatible cohomological classes
satisfying the corresponding condition of
Theorem~\ref{Thm:levelcomp}. To this end, we use
constructions of sect.~\ref{subsect:motivdef}.
For the virtual properties ($\ve$) and ($\vf$)
this step of the proof is more difficult, and it is
discussed in sect.~\ref{subsect:virtual}.

\begin{Pro}\label{Pro:compnonvirt} Any of the properties
{\rm ($\imm$)}, {\rm ($\hi$)}, {\rm ($\ebd$)}, {\rm ($\fib$)}, {\rm ($\npc$)}
of a manifold
$M\in\fM$
implies the existence of a compatible collection of
cohomological classes satisfying the corresponding condition
of Theorem~\ref{Thm:levelcomp}.
\end{Pro}

\begin{proof} For
($\imm$) and ($\npc$)
this follows from Lemmas~\ref{Lem:motivcohom1}
and \ref{Lem:motivcohom2} respectively.

($\hi$) If an immersion
$g:S\to M$
is horizontal, then for the corresponding compatible
cohomological classes
$\set{l_v}{$v\in V$}$
we have
$l_v(f_v)>0$
by the construction.

($\ebd$) Let
$g:S\to M$
be a
$\pi_1$-injective
embedding. Consider the corresponding
compatible collection
$\set{l_v\in H^1(M_v;\R)}{$v\in V$}$
and assume that
$l_w(f_w)\cdot l_{-w}(f_{-w})\neq 0$
for some edge
$w\in W$,
where
$l_w=i_w^\ast l_v$.
Then the embedding
$g$
is horizontal in both blocks
$M_v$, $M_{v'}\sub M$, $v=w^-$, $v'=w^+$.
In this case, the intersection
$g(S)\cap T_{|w|}$
consists of a collection of parallel circles. Hence, they all are
consistent or not simultaneously. Thus
$l_w=\pm l_{-w}$.

($\fib$) Let
$p:M\to S^1$
be a fibration. We fix a fiber
$S=p^{-1}(t)$, $t\in S^1$.
One can assume that the surface
$S\sub M$
is horizontal. As usual, using this surface, we define
compatible cohomological classes
$\set{l_v}{$v\in V$}$.
It follows from horizontality that
$l_v(f_v)>0$
for all
$v\in V$.
We fix a generator
$\al$
of the group
$H^1(S^1;\Z)$
and put
$\ep_v=\sgn\al(p_\ast(f_v))\in\{\pm 1\}$.
Then
$l_v=\ep_vp_v^\ast\al$,
where
$p_v$
is the restriction of
$p$
to the maximal block
$M_v$.
Thus for every edge
$w\in W$
we have
$$l_{-w}(f_w)=i_{-w}^\ast l_{v'}(f_w)
   =i_{-w}^\ast\ep_{v'}p_{v'}^\ast\al(f_v)=
   \ep_{v'}p_v^\ast\al(f_v)=
   \ep_{v'}\ep_vl_w(f_w),$$
where
$v$, $v'$
are the vertices of
$w$.
\end{proof}

\subsubsection{Finding a
$\npc$-metric
on a graph-manifold}\label{subsubsect:constrnpc}

To complete the proof of Theorems~\ref{Thm:levelcomp} and
\ref{Thm:bknlevel}
for the property
($\npc$),
it remains to construct a
$\npc$-metric
on
$M$,
if a compatible collection of cohomological classes
$\set{l_v}{$v\in V$}$
satisfying the condition
($\npc$)
of Theorem~\ref{Thm:levelcomp} is given. We do that using
Lemma~\ref{Lem:locextnpc}.

For
$w\in\d v$
we put
$l_w=i_w^\ast l_v\in H^1(T_{|w|};\R)$
and define a symmetric bilinear form
$g_w$
on
$H_1(T_{|w|};\R)$
by
$g_w(f_w,f_w)=(l_w(f_w))^2$,
$g_w(f_{-w},f_{-w})=(l_{-w}(f_{-w}))^2$,
and
$g_w(f_w,f_{-w})=l_w(f_w)\cdot l_w(f_{-w})$.
It follows from the
($\npc$)-condition
of compatibility
$|l_{-w}(f_w)|<l_w(f_w)$
that the form
$g_w$
is positive definite. Furthermore, it obviously satisfies
the condition (1) of Lemma~\ref{Lem:locextnpc}. The condition (2)
of the same Lemma is also fulfilled since the forms
$g_w$
are defined via the classes
$l_v$.
By Lemma~\ref{Lem:locextnpc}, we obtain a
$\npc$-metric
on every maximal block
$M_v$.
It follow from the
($\npc$)-condition of symmetry
$l_{-w}(f_w)\cdot l_{-w}(f_{-w})=
  l_w(f_w)\cdot l_w(f_{-w})$
that the constructed metrics coincide on the gluing tori
(cf. the proof of Lemma~\ref{Lem:motivcohom2}). Therefore,
we obtain a
$\npc$-metric
on
$M$.
(This metric is
$C^1$-smooth
on
$M$
and real analytic inside of every block
$M_v$
being modelled on
$\hyp^2\times\R$.
Its curvature is nonpositive in Alexandrov' sense.
The metric can easily be smoothed to an
$C^\infty$-smooth
$\npc$-metric
(see, for example, \cite{L}), though.
\qed

\subsubsection{Rational approximation}\label{subsubsect:ratap}

To implement the properties
($\imm$) -- ($\vf$),
one needs to have compatible collections of integer cohomological
classes while we have only real ones. Thus we first
approximate real classes by rational ones, and then we obtain
integer classes from them. It is convenient to
construct an approximation in terms of solutions to the BKN-equation.
The next proposition is sufficient for the properties
($\imm$), ($\hi$), ($\ebd$), ($\fib$).
To implement the virtual properties
($\ve$), ($\vf$)
one needs to have a symmetric approximating angle function
$\ga'$,
and a symmetric approximation is achieved by another method,
see sect.~\ref{subsect:symratap}.

\begin{Pro}\label{Pro:ratapprox} Let
$(a,\ga)$
be a (real) compatible solution to the BKN-equation over
$M\in\fM$.
Then it can be approximated by rational compatible solutions
$(a',\ga')$
either so that

$\bullet$ $\ga'=\ga$
and
$a_v'=a_v$
for the vertices
$v\in V$
for which
$a_v=0$,
if the angle function
$\ga$
takes values in
$\{0,\pm 1\}$;

\noindent
or so that

$\bullet$ $a_v'>0$
for all
$v\in V$
and
$|\ga_w'|<1$
for all edges
$w\in W$,
if
$|\ga_{w_0}|<1$
for some edge
$w_0\in W$.
\end{Pro}

In other words, if the angle function
$\ga$
takes values on
$\{0,\pm 1\}$
and
$\ga_w=0$
for some edge
$w\in W$,
then we have two types of a rational approximation of
the solution
$(a,\ga)$.
The approximation described in the first part will be used for
the properties
($\imm$), ($\hi$), ($\ebd$)
and
($\fib$).
The approximation described in the second part will be used for
the properties
($\imm$)
and
($\hi$).

The first part is easy to prove: unknown lengths
$a_v$, $v\in V$,
when (rational) angles
$\ga$
are fixed, satisfy a linear homogeneous system of
$|V|$
equations with rational coefficients and the zero
determinant, which follows from the BKN-equation.
By a classical theorem from linear algebra, the solutions
to such a system can be parameterized by linear functions
with rational coefficients. This easily implies a required
approximation. The proof of the second part we start with
construction of perturbations of the initial solution
(in the class of solutions) do not caring first about the
rationality. This is achieved in the following three lemmas.

\begin{Lem}\label{Lem:pertub1} Let
$(a,\ga)$
be a compatible solution to the BKN-equation. Then it
can be approximated by compatible solutions
$(a',\ga')$
so that
$a_v'>0$
for all
$v\in V$.
\end{Lem}

\begin{proof} If the length function
$a$
has zeroes on
$V$,
then since the solution is nontrivial there is a vertex
$v\in V$
with
$a_v>0$
and
$a_{v'}=0$
for a neighboring vertex
$v'$.
Consider an edge
$w\in\d v$
coming in
$v'$.
Recall that then
$\ga_{w'}=\ga_{-w'}=0$
for all
$w'\in\d v'$
since the solution
$(a,\ga)$
is compatible. We choose an arbitrarily small
$a_{v'}'>0$
and find
$\ga_{-w}'$
with
$k_{v'}a_{v'}'=\frac{\ga_{-w}'}{|b_{-w}|}a_v$.
Since the angle function
$\ga$
vanishes on all edges going out of the vertex
$v'$,
we achieve by this the BKN-equation over
$v'$
for the perturbed functions
$a'$, $\ga'$.
The BKN-equation for them is also fulfilled over the vertex
$v$
because the function
$\ga$
vanishes on all edges pointing out from
$v$
to
$v'$.
The procedure obviously leads to a required approximation.
\end{proof}

Using this Lemma, we assume from now on that the length
function
$a$
is positive.

\begin{Lem}\label{Lem:pertub2} Let
$(a,\ga)$
be a compatible solution to the BKN-equation. Assume that
for a vertex
$v\in V$
there is an edge
$w_0\in\d v$
with
$|\ga_{w_0}|<1$.
Then
$(a,\ga)$
can be approximated by compatible solutions
$(a',\ga')$
so that
$|\ga_w'|<1$
for all
$w\in\d v$.
\end{Lem}

\begin{proof} Unknown angles
$\ga_w$, $w\in\d v$
enter only the BKN-equation over the vertex
$v$:
$$k_va_v=\sum_{w\in\d v}\frac{\ga_w}{|b_w|}a_{w^+}.$$
Thus changing only them and keeping this equality hold,
we obtain a new solution to the BKN-equation. One needs
to change only those
$\ga_w$, $w\in\d v$,
for which
$|\ga_w|=1$.
Using that the length function
$a$
is positive, their absolute values can
be (arbitrarily small) decreased at the expense of
$\ga_{w_0}$.
\end{proof}

\begin{Lem}\label{Lem:pertub3} Let
$(a,\ga)$
be a compatible solution to the BKN-equation. Assume that
$|\ga_w|<1$
for all
$w\in\d v$
for some vertex
$v\in V$.
Then
$(a,\ga)$
can be approximated by compatible solutions
$(a',\ga')$
so that
$|\ga_{\pm w}'|<1$
for all
$w\in\d v$.
\end{Lem}

\begin{proof} We change the angle function
$\ga$
on every edge
$-w$
with
$|\ga_{-w}|=1$
and
$w\in\d v$,
and the length function
$a$
on the vertex
$v$
so that
$|\ga_{-w}'|\cdot a_v'=|\ga_{-w}|\cdot a_v$,
where the ratio
$\frac{a_v'}{a_v}>1$
is arbitrarily close to 1 (here
$'$
means taking approximating functions).
By that one achieves the property
$|\ga_{-w}'|<1$
for all edges
$w\in\d v$
with the BKN-equation satisfied for all final vertices
of these edges (except the loops). To fulfill the
BKN-equation for the initial vertex
$v$,
we distribute the excess
$k_v(a_v'-a_v)$
over the edges
$w\in\d v$,
changing in an appropriate way the values
$\ga_w$.
One can do this keeping the condition
$|\ga_w|<1$,
because the value
$|k_v(a_v'-a_v)|$
can be chosen arbitrarily small.
\end{proof}

\begin{proof}[Proof of Proposition~\ref{Pro:ratapprox}]
As it was explained above, to prove the first part of
the Proposition, one can use a rational parameterization
of the solutions to the BKN-equation obtained when
rational values of the angle function
$\ga$
are fixed. It follows from compatibility of the solution
$(a,\ga)$
that the equations over the vertices
$v\in V$
with
$a_v=0$
put no restriction on any approximating function
$a'$
having the same zeroes as the function
$a$.
Thus it remains to prove the second part of
the Proposition.

It suffices to approximate the solution
$(a,\ga)$
by rational solutions
$(a',\ga')$,
satisfying the conditions of the Proposition.
Then the solutions
$(a',\ga')$
sufficiently close to
$(a,\ga)$
will be compatible automatically. According Lemmas~\ref{Lem:pertub1},
\ref{Lem:pertub2}, \ref{Lem:pertub3}
one can assume that the initial solution
$(a,\ga)$
to the BKN-equation satisfies the conditions
$a_v>0$
for all vertices
$v\in V$
and
$|\ga_w|<1$
for all edges
$w\in W$.

Now, we approximate the solution
$(a,\ga)$
by solutions
$(a',\ga')$,
where the length function
$a'$
takes rational values. Fix
$v\in V$
and approximate the number
$a_v$
by rational numbers
$a_v'>0$.
The excess
$k_v(a_v'-a_v)$
from the BKN-equation over the vertex
$v$
we distribute between the summands
$\frac{\ga_w}{|b_w|}a_{w^+}$, $w\in\d v$,
changing only if necessarily appropriate
$\ga_w$,
keeping the condition
$|\ga_w|<1$
and leaving invariant the numbers
$a_{w^+}$.
By this, one achieves the BKN-equation to hold over
the vertex
$v$.
For any
$v'$,
neighboring to
$v$,
only summands from the right hand side of the BKN-equation
over
$v'$,
corresponding to the edges
$w\in\d v'$
with
$-w\in\d v$
are changed during the described manipulations. We compensate
these changes by the condition
$\ga_w'\cdot a_v'=\ga_w\cdot a_v$
for every such edge
$w$.
We have
$|\ga_w'|<1$,
if the ratio
$\frac{a_v}{a_v'}$
is sufficiently close to 1. Therefore, one can assume
that the length function
$a$
from the initial compatible solution
$(a,\ga)$
takes positive rational values on all vertices
$v\in V$.
Now, it is easy to approximate the numbers
$\ga_w$
by rational numbers
$\ga_w'$
for any edge
$w\in W$
fulfilling the BKN-equation.
\end{proof}

\subsubsection{Finding a horizontal immersion}
\label{subsubsect:constrfib}

To complete the proof of Theorems~\ref{Thm:levelcomp} and
\ref{Thm:bknlevel} for the properties
($\imm$)
and
($\hi$),
we construct an immersed horizontal surface in
$M$,
if there is a compatible collection of cohomological
classes on
$M$.
Using Proposition~\ref{Pro:ratapprox}, we find
a compatible rational solution
$(a,\ga)$
to the BKN-equation for which
$a_v>0$
for all
$v\in V$,
and if
$|\ga_{w_0}|=1$
for some edge
$w_0\in W$,
then
$\ga_{-w}=\ga_w=\pm 1$
for all edges
$w\in W$.
This gives an integer collection
$\set{l_v\in H^1(M_v;\Z)}{$v\in V$}$
of compatible classes, which satisfies the following,
more strong, condition of compatibility:
if
$|l_{-w_0}(f_{w_0})|=l_{w_0}(f_{w_0})$
for some edge
$w_0\in W$,
then
$l_{-w}=\pm l_w$
for all edges
$w\in W$.

For every edge
$w\in W$
we put
$l_w^+=l_w+l_{-w}$, $l_w^-=l_w-l_{-w}$
and consider the elements
$c_w^+$, $c_w^-\in H_1(T_{|w|};\Z)$
dual to
$l_w^+$, $l_w^-$.
These elements satisfy the conditions of
Lemma~\ref{Lem:locextimm} because:

\noindent
$\bullet$ $\al_w^{\pm}=c_w^{\pm}\land_wf_w=
         l_w(f_w)\pm l_{-w}(f_w)\ge 0$;

\noindent
$\bullet$
the class
$c_w=c_w^++c_w^-$
is even being dual to the class
$2l_w$;

\noindent
$\bullet$ $c_w\land_wf_w=2l_v(f_v)>0$
for all
$v\in V$
and
$w\in\d v$;
the condition (2) of Lemma~\ref{Lem:locextimm} is also
holds since the classes
$l_w$, $w\in\d v$
are boundary values for the class
$l_v$, $i_w^\ast l_v=l_w$,
see Lemma~\ref{Lem:extension}.

Moreover, the condition of the strong compatibility
means that either
$\al_w^+\cdot\al_w^-\neq 0$
for all
$w\in W$
or
$\al_w^+\cdot\al_w^-=0$
for all
$w\in W$.
Since
$c_{-w}^+=c_w^+$, $c_{-w}^-=-c_w^-$
and if
$\al_w^+\cdot\al_w^-=0$,
then
$c_{-w}=\pm c_w$,
horizontal immersions
$g_v:S_v\to M_v$
provided by Lemma~\ref{Lem:locextimm} for each
maximal block
$M_v$
can be easily pasted on every gluing torus
$T_{|w|}$
at the expense of an appropriate multiplicity
$d$
universal for all blocks. This gives a required
horizontal immersion
$g:S\to M$.
\qed

\begin{Cor}\label{Cor:imm_hi} The property
{\rm ($\imm$)}
implies the property
{\rm ($\hi$).}
\qed
\end{Cor}

\subsubsection{Finding a fibering over the circle}
\label{subsubsect:fib}

To complete the proof of Theorems~\ref{Thm:levelcomp} and
\ref{Thm:bknlevel} for the property
($\fib$),
we construct a fibration
$M\to S^1$
with fiber which is a closed surface of negative Euler
characteristic, if there are compatible cohomological
classes
$\set{l_v}{$v\in V$}$
on
$M$
satisfying the condition
($\fib$)
of Theorem~\ref{Thm:levelcomp}. This condition implies that
$\ep_vl_w=\ep_{v'}l_{-w}$
for every edge
$w\in W$
from
$v$
to
$v'$.
According to the first part of Proposition~\ref{Pro:ratapprox},
one can assume the collection
$\set{l_v}{$v\in V$}$
to be integer. Using the canonical isomorphism
$H^1(M_v;\Z)=[M_v,S^1]$
and that
$l_v(f_v)\neq 0$,
we find for every maximal block
$M_v$
a map
$M_v\to S^1$
corresponding to the class
$\ep_vl_v$
of nonzero degree on the Seifert fibers for
$M_v$.
It is well known that such a map is homotopic to
a fibration
$p_v:M_v\to S^1$
(see the discussion of this question in \cite{N2}).
Since
$\ep_vl_w=\ep_{v'}l_{-w}$
for any edge
$w\in W$
from
$v$
to
$v'$,
we can choose the fibrations
$p_v$
and
$p_{v'}$
coinciding on the gluing torus
$T_{|w|}$.
This gives a required fibration
$p:M\to S^1$.

\subsubsection{Finding an embedded surface}\label{subsubsect:ebd}

We consider now the implementation of the property
($\ebd$).
By the first part of Proposition~\ref{Pro:ratapprox},
one can assume that all classes
$l_v\in H^1(M_v;\Z)$, $v\in V$
of the corresponding compatible collection are even.
Given a block
$M_v$
with
$l_v(f_v)\neq 0$,
we find as above a fibration
$M_v\to S^1$.
Its fibers are horizontally embedded surfaces whose
relative class corresponds to the class
$l_v$
under the isomorphism
$H_2(M_v,\d M_v;\Z)\simeq H^1(M_v;\Z)$.
If
$l_v(f_v)\cdot l_{v'}(f_{v'})\neq 0$
for neighboring vertices
$v$, $v'$,
then
$l_w=\pm l_{-w}$
by the condition
($\ebd$)
of Theorem~\ref{Thm:levelcomp}, where the edge
$w\in W$
points out of
$v$
to
$v'$.
Thus one can glue corresponding horizontal surfaces in the blocks
$M_v$, $M_{v'}$
on the separating torus
$T_{|w|}$
(slightly deforming them if necessarily). This gives
a surface
$S'$
horizontally embedded into the union of all
maximal blocks
$M_v$
with
$l_v(f_v)\neq 0$.

Assume now that
$l_v(f_v)=0$
for some vertex
$v\in V$.
We can further assume that
$l_{v'}(f_{v'})\neq 0$
for a neighboring vertex
$v'$.
Consider an edge
$w\in W$
from
$v$
to
$v'$.
The compatibility condition
$|l_{-w}(f_w)|\le l_w(f_w)=0$
implies that the horizontal surface
$S'\cap M_{v'}$
meets the torus
$T_{|w|}$
by a collection of circles parallel to the Seifert
fibers of
$M_v$.
It follows from the parity condition that the number of
the circles is even. It is not difficult to close these
circles by vertically embedded annuli in
$M_v$
so that the annuli are not parallel to the block boundary
and pairwise disjoint (taking into account all edges
$w\in\d v$).
We leave details to the reader. The described procedure
obviously gives a closed surface embedded into
$M$,
whose intersection with any block either is
horizontal or consists of vertical annuli, or is empty.
It is known (see \cite[Lemma~3.4]{RW}) that any such
an embedding is
$\pi_1$-injective.
\qed

This completes the proof of Theorems~\ref{Thm:levelcomp}
and \ref{Thm:bknlevel} for all properties except the virtual
ones,
($\ve$)
and
($\vf$).


\subsection{Virtual properties: necessary conditions}\label{subsect:virtual}

We start the proof of Theorems~\ref{Thm:levelcomp},
\ref{Thm:bknlevel} for the virtual properties
($\ve$)
and
($\vf$).
The proof of necessity of the corresponding conditions is
based on the existence of characteristic coverings of the
graph-manifolds and Descent-Lemma~\ref{Lem:down}. This Lemma
also implies the existence of a
$\npc$-metric
on a graph-manifold if some its finite cover carries a
$\npc$-metric.
This fact has been proved earlier by another method in \cite{KL}
for a wider class of manifolds.

\subsubsection{Descent Lemma}\label{subsubsect:down}

Let
$s\ge 1$
be an integer. A covering
$p:\wt T\to T$
of tori is called
$s$-characteristic,
if the image of the induced homomorphism
is
$p_\ast(H_1(\wt T;\Z))=s\cdot H_1(T;\Z)$.
A covering
$\wt M\to M$
of manifolds of the class
$\fM$
is called
$s$-characteristic,
if its restriction to each JSJ-torus in
$\wt M$
is
$s$-characteristic.

The existence of
$s$-characteristic
coverings of a manifold
$M\in\fM$
with an arbitrary multiple
$s$
of some
$s_0$
has been proved in \cite{LW} (by arguments in spirit of
the proof of Lemma~\ref{Lem:locextimm}), see also \cite{N2}.
Furthermore, the covering can be pushed through any given
finite covering
$\wh M\to M$.

Characteristic covering are good by that they allow easily
to trace the behavior of the intersection indices and the
charges. Namely, let
$\wt M\to M$
be a
$s$-characteristic
covering of manifolds of the class
$\fM$.
Denote by
$p:\wt\Ga\to\Ga$
the induced map of the graphs
$\wt\Ga=\Ga_{\wt M}$, $\Ga=\Ga_M$.
The restriction
$\wt M_{\wt v}\to M_v$
of the covering to any maximal block
$\wt M_{\wt v}\sub\wt M$
has the degree
$d_{\wt v}\cdot s^2$,
the number
$s^2\cdot\sum_{p(\wt v)=v}d_{\wt v}$
is independent of the vertex
$v\in V$
of the graph
$\Ga$
and it equals the degree of the covering
$\wt M\to M$.

We assume that an orientation of
$M$
and fiber orientations of its maximal blocks are fixed,
and for
$\wt M$
the orientations induced by the covering are chosen
(it can be done because
$M$
contains no one-sided Klein bottles). It follows from
the characteristic property that the intersection indices
for the edges
$\wt w\in\wt W$
of
$\wt\Ga$
and
$w=p(\wt w)\in W$
of
$\Ga$
coincide,
$b_{\wt w}=b_w$.
For the charges of
$\wt v\in\wt V$
and
$v=p(\wt v)\in V$
one holds
$$k_{\wt v}=d_{\wt v}\cdot k_v,$$
what easy to see from the fact that every boundary torus
$T_w$
of the block
$M_v$
is
$s$-characteristically
covered by
$d_{\wt v}$
tori of the block
$\wt M_{\wt v}$,
and thus their contribution into the charge
$k_{\wt v}$
is
$d_{\wt v}$-fold
of the contribution of
$T_w$
into the charge
$k_v$
(one can find details in \cite{LW}, \cite{N2}).
These two properties immediately imply the following Lemma.

\begin{Lem}\label{Lem:up0} Let
$p:\wt\Ga\to\Ga$
be the map of the graphs
$\wt\Ga=\Ga_{\wt M}$, $\Ga=\Ga_M$,
of manifolds
$M$, $\wt M\in\fM$
induced by a characteristic covering
$\wt M\to M$.
Then for any solution
$(a,\ga)$
to the BKN-equation over
$M$
the functions
$\wt a=a\circ p$, $\wt\ga=\ga\circ p$
form a solution to the BKN-equation over
$\wt M$.
\qed
\end{Lem}

Furthermore, only these two properties of
characteristic coverings will be used in the proof of the
following Descent Lemma.

\begin{Lem}\label{Lem:down} Let
$\wt M\to M$
be a characteristic covering of a manifold
$M\in\fM$.
If the BKN-equation over
$\wt M$
has a compatible solution satisfying one of the conditions
{\rm ($\ve$)}, {\rm ($\vf$)}, {\rm ($\npc$)} of
Theorem~\ref{Thm:bknlevel}, then the BKN-equation over
$M$
has a compatible solution satisfying the same condition.
\end{Lem}

\begin{proof} We use notations introduced above for
a characteristic covering, and for every vertex
$v\in V$
(of the graph of
$M$)
we define a quadratic form
$D_v$
on
$\R^{\wt V}$,
$$(D_vx,x)=\sum_{p(\wt v)=v}d_{\wt v}\cdot x_{\wt v}^2,
  \quad x\in\R^{\wt V}.$$
Let
$(\wt a,\wt\ga)$
be a compatible solution to the BKN-equation over the cover
$\wt M$.
For every edge
$w\in W$
(of the graph of
$M$)
we define a quadratic form
$G_w$
on
$\R^{\wt V}$,
putting
$$(G_wx,x)=\sum_{p(\wt w)=w}\wt\ga_{\wt w}x_{\wt w^-}x_{\wt w^+},
  \quad x\in\R^{\wt V}.$$
From the condition
$|\wt\ga_{\wt w}|\le 1$,
we obtain
$(G_wx,x)^2\le(D_vx,x)(D_{v'}x,x)$
for every
$x\in\R^{\wt V}$,
where
$v=w^-$
is the initial vertex of the edge
$w$, $v'=w^+$
its final vertex. To deduce the last inequality, one needs
to apply twice the Cauchy-Schwartz inequality and to use that
the number of edges
$\wt w\in\d\wt v$
projected onto an edge
$w$, $p(\wt w)=w$,
is equal to
$d_{\wt v}$
for every vertex
$\wt v\in\wt V$
with
$p(\wt v)=v$.

Using the relations for the intersection indices and
the charges between
$M$
and its cover
$\wt M$,
and the BKN-equation
$$k_{\wt v}\wt a_{\wt v}=\sum_{\wt w\in\d\wt v}
   \frac{\wt\ga_{\wt w}}{|b_{\wt w}|}\wt a_{\wt w^+},\quad
    \wt v\in\wt V,$$
for
$(\wt a,\wt\ga)$,
we obtain
$$k_v(D_v\wt a,\wt a)=\sum_{w\in\d v}\frac{1}{|b_w|}(G_w\wt a,\wt a),
  \quad v\in V.$$
Thus the functions
$a:V\to\R$, $\ga:W\to\R$,
defined by
$a_v=(D_v\wt a,\wt a)^{1/2}$
(the form
$D_v$
is nonnegative),
$$\ga_w=\frac{(G_w\wt a,\wt a)}{a_{w^-}\cdot a_{w^+}},\quad
\text{if}\quad a_{w^-}\cdot a_{w^+}\neq 0,$$
and
$\ga_w=0$,
if
$a_{w^-}\cdot a_{w^+}=0$,
make up a compatible solution to the BKN-equation for
$M$: $|\ga_w|\le 1$
for all
$w\in W$.

To complete the proof, we note that
the length function
$a$
is positive, if
$\wt a$
is, and the angle function
$\ga$
is symmetric, if
$\wt\ga$
is
($\wt\ga_{-\wt w}=\wt\ga_{\wt w}$
for all
$\wt w\in\wt W$).
Furthermore, if
$|\wt\ga|<1$,
then
$|\ga|<1$.
Hence, each of the properties
($\ve$), ($\vf$), ($\npc$)
from Theorem~\ref{Thm:bknlevel} of the BKN-equation over
$\wt M$
implies the corresponding property of the BKN-equation over
$M$.
\end{proof}

\subsubsection{Existence of a solution to BKN-equation}

\begin{Pro}\label{Pro:bknvirt} Each of the properties
{\rm ($\ve$)}, {\rm ($\vf$)}
of a manifold
$M\in\fM$
implies the existence of a compatible solution to
the BKN-equation over
$M$,
which satisfies the corresponding condition of
Theorem~\ref{Thm:bknlevel}.
\end{Pro}

\begin{proof} Assume that a manifold
$M\in\fM$
has one of the properties
($\ve$), ($\vf$).
There is a finite covering
$\wt M\to M$,
such that the manifold
$\wt M$
has the corresponding property
($\ebd$)
or
($\fib$).
One can suppose that the covering is characteristic. Then
by Proposition~\ref{Pro:compnonvirt} and due to the fact that
Theorems~\ref{Thm:levelcomp}, \ref{Thm:bknlevel}
are equivalent, the BKN-equation over
$\wt M$
has a compatible solution
$(\wt a,\wt\ga)$,
which satisfies the corresponding property
($\ebd$)
of
($\fib$)
from Theorem~\ref{Thm:bknlevel}. In any case, the angle
function
$\wt\ga$
is symmetric, and in the
($\fib$)
case the length function
$\wt a$
is positive. Now, Lemma~\ref{Lem:down} provides the
existence of a required solution to the BKN-equation over
$M$.
\end{proof}

\subsection{Implementing the virtual properties}
\label{subsect:sufficevirtprop}

Here we prove that the conditions of Theorem~\ref{Thm:bknlevel}
are sufficient to implement the virtual properties.
We start with the key Ascent Lemma.

\subsubsection{Ascent Lemma}\label{subsubsect:up}

Let
$g:S\to M$
be a horizontal immersion of a closed surface
$S$
into a manifold
$M$
of the class
$\fM$, $\cT$
be the JSJ-surface in
$M$.
Its preimage
$\cT_g=g^{-1}(\cT)$
consists of a finite, disjoint collection of simple,
closed, noncontractible curves on
$S$.
One can assume that every curve
$T\sub\cT_g$
has as the image a simple, closed curve
$g(T)$
on the corresponding torus from
$\cT$.

Let
$\Ga_g$
be the graph dual to the decomposition of the surface
$S$, $S=\cup_{v\in V_g}S_v$,
defined by
$\cT_g$.
In other words, the vertex set
$V_g$
of the graph
$\Ga_g$
is the set of the connected components of the splitting
$S|\cT_g$
(which are called the blocks
$S_v$
of
$S$),
and the set of (nonoriented) edges can be identified with
the collection
$\cT_g$.
The set of oriented edges of the graph
$\Ga_g$
is denoted by
$W_g$,
and for the edges
$w\in W_g$
we shall use the same agreements and notations as for
the oriented edges
$w\in W$
of the graph
$\Ga=\Ga_M$.
The immersion
$g$
induces the map
$\ov g:\Ga_g\to\Ga$
of the graphs taking the vertices into the vertices and
edges into the edges.

Fix an orientation of
$M$,
fiber orientations of its maximal blocks, orientations of
the curves
$T\sub\cT_g$
and assign to each edge
$w\in W_g$
the number
$d_w=|[g(T)]\land_{\ov g(w)}f_{\ov g(w)}|$,
where the curve
$T\sub\cT_g$
corresponds to the nonoriented edge
$(w,-w)$
of the graph
$\Ga_g$.
In other words,
$d_w$
is equal to the degree with which the curve
$T$
is mapped onto the corresponding boundary component
of the base orbifold
$\cO M_{v'}$
of the block
$M_{v'}\sub M$,
into which the block
$S_{w^-}$
of the surface
$S$,
corresponding to the initial vertex of
$w$,
is mapped under the composition
$\pi_{v'}\circ g$,
where
$\pi_{v'}:M_{v'}\to\cO M_{v'}$
is the canonical projection.

All oriented cycles of the graph
$\Ga_g$
we consider below consist of oriented edges of that graph
(with compatible orientations).

\begin{Lem}\label{Lem:up} The immersion
$g$
is a virtual embedding, i.e.,
$g=\wt g\circ p$
for some finite covering
$p:\wt M\to M$
and a horizontal embedding
$\wt g:S\to\wt M$,
if and only if the following condition of
cyclic balance holds:
$$\prod_{w\in c}\frac{d_w}{d_{-w}}=1$$
for every oriented cycle
$c$
in the graph
$\Ga_g$.
\end{Lem}

\begin{proof}[Sketch of the proof] It is easy to see, using that
the maximal blocks of
$M$
are Seifert bundles, that the cover space of the covering
associated with subgroup
$g_\ast(\pi_1(S))\sub\pi_1(M)$
is homeomorphic to the product
$S\times\R$.
One can further assume that the covering
$q:S\times\R\to M$
itself coincides over the zero section
$S\times\{0\}\sub S\times\R$
with the immersion
$g$.
In what follows, we identify
$S$
with
$S\times\{0\}$.

For every curve
$T\sub\cT_g$,
the infinite cylinder
$T\times\R\sub S\times\R$
covers the corresponding torus
$T_{|\ov g(w)|}\sub\cT$,
where
$T=(w,-w)$.
Since
$g(T)\sub T_{|\ov g(w)|}$
is a simple closed curve, the preimage
$q^{-1}(g(T))\cap T\times\R$
consists of a countable collection of parallel copies of
the curve
$T$.
These copies can be labeled by integers
$n\in A\sub\Z$
so that any finite subcylinder in
$T\times\R$,
bounded by
$T$
and its copy
$T_{w,n}$,
is mapped into the torus
$T_{|\ov g(w)|}$
with degree
$|n|\cdot d_w$.
We shall say that
$T_{w,n}$
is the copy of the level
$n$.
Then
$T=T_{w,0}$.

The curve
$T$
is a boundary component of a block
$S_v\sub S$,
where
$w\in\d v\sub V_g$.
Then its parallel copy
$T_{w,n}$
is a boundary component of a copy
$S_v'\sub S_v\times\R$
of the block
$S_v$,
and
$q(S_v')=g(S_v)$.
Using the Seifert bundle structure over the block
$M_{\ov g(v)}\sub M$
and that the degree is invariant under a homotopy, it is
easy to see that the level of every component of the boundary
$\d S_v'$
is one and the same. Thus we have a collection
$S_{v,n}$, $n\in A$,
of copies of the block
$S_v$,
which lie in the cylinder
$S_v\times\R$
with images
$q(S_{v,n})=g(S_v)$
for every vertex
$v\in V_g$.
Furthermore, for each adjacent vertices
$v=w^-$, $v'=w^+$
and copies
$S_{v,n}$, $S_{v',n'}$
of the corresponding blocks
$S_v$, $S_{v'}$,
adjacent along a common boundary component,
we have a key relation
\[\label{eq:degree}
n\cdot d_w=n'\cdot d_{-w},\tag{\rm deg}\]
in which the left and the right hand sides are different
representations of the degree of the map of the corresponding
cylinder in the torus
$T_{|\ov g(w)|}$.

If
$g$
is a virtual embedding, then the manifold
$\wt M$
from the corresponding covering
$\wt M\to M$,
which fibers over the circle, can be covered in turn by a product
$S\times S^1$,
because the holonomy of the gluing of the manifold
$\wt M$
from the cylinder
$S\times[0,1]$
has a finite order (up to an isotopy).
Thus the covering
$q:S\times\R\to M$
can be pushed through the covering
$S\times S^1\to\R$,
and hence the preimage
$q^{-1}(g(S))$
contains a parallel copy
$S'$
of the surface
$S$.
Its block decomposition is
$S'=\cup_{v\in V_g}S_{v,n(v)}$,
and according what was said above, for every edge
$w\in W_g$
we have
$n_{w^-}\cdot d_w=n_{w^+}\cdot d_{-w}$.
Thus for every oriented cycle
$c\sub\Ga_g$
we have
$$\prod_{w\in c}\frac{d_w}{d_{-w}}=
  \prod_{w\in c}\frac{n_{w^+}}{n_{w^-}}=1,$$
since
$w^+=w'^-$
for each consecutive edges
$w$, $w'$
of the cycle
$c$.

Conversely, suppose that for a horizontal immersion
$g:S\to M$
the cyclic balance holds. Together with the relation
(\ref{eq:degree}),
it allows to find, starting with some copy
$S_{v,n}$, $n\in A\sm\{0\}$,
and using the obvious continuation method, a copy of the
surface
$S$
in the preimage
$q^{-1}(g(S))\sub S\times\R$
distinct from
$S$.
This implies the existence of a required covering
$p:\wt M\to M$
and of an embedding
$\wt g:S\to\wt M$
with
$\wt g\circ p=g$.
\end{proof}

\subsubsection{Compatible symmetric solutions}
\label{subsubsect:compsymsol}

Here we formulate a proposition which will be used
to complete the proof of Theorems~\ref{Thm:levelcomp} and
\ref{Thm:bknlevel}.

A compatible symmetric solution
$(a,\ga)$
to the BKN-equation over a manifold
$M\in\fM$
is said to be a NPC-solution, if it satisfies the condition
($\npc$)
of Theorem~\ref{Thm:bknlevel}, i.e., the length function
is positive
$a>0$,
and the angle function
$|\ga|<1$.

\begin{Pro}\label{Pro:symcompsol} If the BKN-equation over
$M\in\fM$
has a compatible, symmetric solution
$(a,\ga)$,
then the same equation has a compatible, symmetric solution
$(a',\ga')$
which satisfies one of the following two conditions:

$(1)$ $(a',\ga')$
is a NPC-solution;

$(2)$
the angle function
$\ga'$
takes values in
$\{0,\pm 1\}$,
and
$|\ga_w'|=1$
for every edge
$w\in W$
with
$a_{w^-}'\cdot a_{w^+}'\neq 0$.

Furthermore, if the length function is positive,
$a>0$,
then the equation has a compatible, symmetric solution
$(a',\ga')$
satisfying (1) or

$(3)$
the length function is positive,
$a'>0$,
and the angle function
$\ga'$
takes values in
$\{0,\pm 1\}$.
\end{Pro}

We prove this Proposition in sect.~\ref{subsect:idealdecomp}.

\begin{Rem}\label{Rem:ratsymcomsol} To implement the virtual
properties
($\ve$)
and
($\vf$)
we need rational solutions to the BKN-equation.
In the cases (2) and (3), the existence of a rational
solution
$(a',\ga')$
follows from the first part of Proposition~\ref{Pro:ratapprox}.
Therefore, the problem of rationality is reduced to the case of
a NPC-solution. There are manifolds
$M\in\fM$,
for which the BKN-equation has a NPC-solution, but has no
rational NPC-solution. However as it is shown in
sect.~\ref{subsect:symratap} (see Proposition~\ref{Pro:symratappr}),
there always is a characteristic cover of such
a manifold for which there is a rational NPC-solution. This suffices
for our purposes because we consider virtual properties
(see Remark~\ref{Rem:irrdipole}).
\end{Rem}

\subsubsection{Implementing
$(\vf)$}\label{subsubsect:vf}

Here we complete the proof of Theorem~\ref{Thm:bknlevel}($\vf$).
Assume that the BKN-equation over a manifold
$M\in\fM$
has a compatible symmetric solution
$(a,\ga)$
with positive length function,
$a>0$.
We show that
$M$
is virtually fibered over the circle.

According to Proposition~\ref{Pro:symcompsol} and
Remark~\ref{Rem:ratsymcomsol}, one can assume, taking
if necessarily an appropriate characteristic covering, that
$(a,\ga)$
is a rational solution. Then it follows from
Proposition~\ref{Pro:compequivbkn} that there are compatible
integer cohomological classes
$\set{l_v}{$v\in V$}$
on
$M$
such that
$l_v(f_v)>0$
for all vertices
$v\in V$
and
$l_{-w}(f_w)\cdot l_{-w}(f_{-w})=
 l_w(f_w)\cdot l_w(f_{-w})$
for every edge
$w\in W$
(we suppose that fiber orientations of maximal blocks and
an orientation of
$M$
are fixed). For edges
$w\in W$
with
$|\ga_w|<1$
we have
$|l_{-w}(f_w)|<l_w(f_w)$,
and for edges with
$|\ga_w|=1$
we have
$l_{-w}=\pm l_w$
due to the symmetry property of
$\ga$
and compatibility of
$\set{l_v}{$v\in V$}$.

Using this collection, we find a horizontal immersion
$g:S\to M$
as in sect.~\ref{subsubsect:constrfib}. Namely, for every edge
$w\in W$
we put
$l_w^+=l_w+l_{-w}$, $l_w^-=l_w-l_{-w}$.
Then the elements
$c_w^+$, $c_w^-\in H_1(T_{|w|};\Z)$,
dual to
$l_w^+$, $l_w^-$,
satisfy the conditions of Lemma~\ref{Lem:locextimm}. Note that
the numbers
$\al_w^{\pm}=c_w^{\pm}\land_wf_w=l_w(f_w)\pm l_{-w}(f_w)$
are positive for edges
$w\in W$
with
$|\ga_w|<1$,
and for edges with
$|\ga_w|=1$
we have
$\al_w^-\cdot\al_w^+=0$.
Furthermore,
$\al_w^-\cdot\al_w^+=0$
if and only if
$\al_{-w}^-\cdot\al_{-w}^+=0$
by symmetry of
$\ga$.

Since
$c_{-w}^+=c_w^+$, $c_{-w}^-=-c_w^-$,
and if
$\al_w^+\cdot\al_w^-=0$,
then
$c_{-w}=\pm c_w$,
the horizontal immersions
$g_v:S_v\to M_v$
provided by Lemma~\ref{Lem:locextimm} for every maximal
block
$M_v$
can be easily pasted together on every gluing torus
$T_{|w|}$
at the expense of an appropriate multiplicity
$d$
universal for all blocks, and they give a required
horizontal immersion
$g:S\to M$.

It follows from the construction of Lemma~\ref{Lem:locextimm}
that the graph
$\Ga_g$,
defined for the immersion
$g$
in sect.~\ref{subsubsect:up}, has the same vertex set as
the graph
$\Ga=\Ga_M$, $V_g=V$,
and its edge set is twice of that of
$\Ga$, $W_g=2W$,
in the sense that to every edge
$w\in W$
it correspond two edges
$w_{\pm}\in W_g$
with the same end points as
$w$.
The number
$d_{w_\pm}$,
defined in sect.~\ref{subsubsect:up}, is equal to
$d\cdot\al_w^{\pm}=d\cdot(l_w(f_w)\pm l_{-w}(f_w))$
in the case
$\al_w^+\cdot\al_w^-\neq 0$,
and
$d_{w_-}=d_{w_+}=d\cdot l_w(f_w)$
in the case
$\al_w^+\cdot\al_w^-=0$,
where the factor
$d>0$
is one and the same for all
$w\in W$.
It follows from the symmetry condition
$l_{-w}(f_w)\cdot l_{-w}(f_{-w})=
 l_w(f_w)\cdot l_w(f_{-w})$,
$w\in W$,
that
$$\frac{d_w}{d_{-w}}=\frac{l_v(f_v)}{l_{v'}(f_{v'})}$$
for all
$w\in W_g$,
where
$v$
is the beginning and
$v'$
is the end of
$w$.

This immediately implies the cyclic balance condition for
$g$.
Hence, by Lemma~\ref{Lem:up},
$g$
is a virtual horizontal embedding, and thus the manifold
$M$
virtually fibers over the circle.
\qed

\subsubsection{Implementing
$(\ve)=(\npc)\cup(\ebd)$}\label{subsubsect:npcebd}

Here we complete the proof of Theorem~\ref{Thm:bknlevel}($\ve$).
Assume that the BKN-equation over a manifold
$M\in\fM$
has a compatible symmetric solution
$(a,\ga)$.
We show that
$M$
has the property
$(\ve)$,
and
$(\ve)=(\npc)\cup(\ebd)$.
According to Proposition~\ref{Pro:symcompsol}
one can assume that the solution
$(a,\ga)$
satisfies one of the conditions (1), (2). In the first case,
taking if necessarily an appropriate characteristic covering,
we suppose also by Remark~\ref{Rem:ratsymcomsol}
that the solution is rational. In the second case, one can
also suppose by the first part of Proposition~\ref{Pro:ratapprox}
that the solution is rational (there is no need to take any
finite covering in this case, what is important for the proof
of the equality
$(\ve)=(\npc)\cup(\ebd)$).

In the first case,
$(a,\ga)$
is a NPC-solution. By sect.~\ref{subsubsect:vf}, the manifold
$M$
virtually fibers over the circle, in particular, it has the property
($\ve$).
In the second case, the solution
$(a,\ga)$
satisfies the condition
($\ebd$)
of Theorem~\ref{Thm:bknlevel}. By sect.~\ref{subsubsect:ebd},
the manifold
$M$
has the property
($\ebd$)
and hence the property
($\ve$).
This also proves the equality
$(\ve)=(\npc)\cup(\ebd)$.
\qed

\subsubsection{Historical remarks} Ascent Lemma~\ref{Lem:up}
is Theorem~2.3 from \cite{RW}, where using this result the authors
give an example of a horizontal immersion of a closed surface
in a graph-manifold which is not a virtual embedding.
It was discovered in \cite{Sv1} that the symmetry condition
of a solution to the BKN-equation implies the cyclic balance
condition for the corresponding horizontal immersion
(sect.~\ref{subsubsect:vf}), and hence the Ascent Lemma can
be applied for the proof that a graph-manifold virtually
fibers over the circle.

One and the same symmetry condition of
solutions to the BKN-equation provides, on the one hand,
that the gluings of maximal blocks are isometric while
constructing a NPC-metric on a graph-manifold, and, on the other
hand, the cyclic balance condition while finding a virtual
embedding of a surface in a graph-manifold. This fact looks somewhat
mystically.

\subsection{Ideals}\label{subsect:idealdecomp}

Here we prove Proposition~\ref{Pro:symcompsol} using
global methods. The notion of an ideal, the decomposition
principle and the polarized discharge Lemma are keys for that.
These notions and results play an important role also for
the proofs of spectral criteria in sect.~\ref{sect:proofspcr}.

A labeled graph
$X=(\Ga,|B|,K)$
is called an ideal if the BKN-equation over every connected
labeled graph
$X'$
containing
$X$
as a subgraph has a NPC-solution. We describe the most important
ideals in the following two lemmas.

\begin{Lem}\label{Lem:posideal} Assume that all charges of
a labeled graph
$X$
have one and the same sign
$\neq 0$.
If
$$(A_X^+a,a):=\sum_{v\in V}|k_v|a_v^2
  -\sum_{w\in W}\frac{a_{w^-}a_{w^+}}{|b_w|}<0$$
for some nonnegative function
$a:V\to\R$,
then
$X$
is an ideal.
\end{Lem}

\begin{Lem}\label{Lem:idipole} Let
$X$
be a labeled graph over the circle with odd number of
edges
$p\ge 1$,
or over a linear graph
$\Ga$
with
$p\ge 2$
edges. Assume that in the circle case at most one vertex
has a nonzero charge, and if
$p=1$,
i.e., in the case of a loop, we require that
$|k\cdot b|<2$,
where
$k$
is the charge, and
$b$
the intersection index. In the case of a linear graph,
we suppose that only extreme vertices have nonzero
charges, which are of the same sign if
$p$
is odd, and of different signs, if
$p$
is even. Then
$X$
is an ideal.
\end{Lem}

We prove these Lemmas in sect.~\ref{subsect:decomprin}. Here,
using them, we prove the following sufficient condition for
the existence of NPC-solutions. Recall that the operator
$H_M:\R^V\to\R^V$
is defined in sect.~\ref{subsubsect:signcomp}.

\begin{Lem}\label{Lem:negeig} Let
$X$
be the labeled graph of an oriented manifold
$M\in\fM$.
If the operator
$H_M$
has a negative eigenvalue or the function
$s:U\to\{0,\pm 1\}$
from its definition is zero, then the BKN-equation over
$X$
has a NPC-solution.
\end{Lem}

\begin{proof} If the sign components graph of
$X$
is not bipartite, then the graph
$X$
contains one of the ideals described in Lemma~\ref{Lem:idipole}.
Thus we suppose that the sign components graph is bipartite.

If the charges
$k_v=0$
for all vertices
$v\in V$,
then the functions
$a\equiv 1$, $\ga\equiv 0$
form a NPC-solution
$(a,\ga)$
of the BKN-equation. Thus we also suppose that the function
$s:U\to\{0,\pm 1\}$
takes values in
$\{\pm 1\}$.

Since the operator
$H_M$
has a negative eigenvalue, there is a sign component
$u\in U$,
for which the operator
$H_u=D_u-J_u$
has a negative eigenvalue. If
$u\in U_0$,
i.e., if
$u$
is a vertex with zero charge, then all edges of the graph
$\Ga_u$
are loops with common vertex
$u$.
In this case, the graph
$X$
contains an ideal by Lemma~\ref{Lem:idipole}. Thus we assume
that
$u\not\in U_0$.
If
$s(u)\sgn(u)=-1$,
where
$\sgn(u)$
is the sign of the charges of the vertices from
$u$,
then the graph
$X$
contains one of the ideals described in Lemma~\ref{Lem:idipole},
because there is a component
$u_0$
with
$s(u_0)\sgn(u_0)=1$.
The proof is completed in this case. Thus one can assume that
$s(u)\sgn(u)=1$.
In this case, the operator
$H_u$
coincides with operator
$A_u^+$,
and the component
$u$
is an ideal by Lemma~\ref{Lem:posideal}.
\end{proof}

\begin{proof}[Proof of Proposition~\ref{Pro:symcompsol}]
Consider a quadratic form
$F:\R^V\to\R$, $F(x)=(H_Mx,x)$.
By Lemma~\ref{Lem:negeig}, one can assume that
$F$
is nonnegative on
$\R^V$,
and the function
$s:U\to\{0,\pm 1\}$
from the definition of the operator
$H_M$
takes values in
$\{\pm 1\}$.

On the other hand, for a compatible symmetric solution
$(a,\ga)$
we have
$F(a)\le 0$
by Lemma~\ref{Lem:symnonposeigen}.
Hence
$F(a)=0$
and
$s(u)\ga_w=1$
for all sign components
$u\in U$
and all edges
$w\in W_u$,
for which
$a_{w^-}\cdot a_{w^+}\neq 0$.
Moreover, the derivative
$\d_vF(a)=0$
for all
$v\in V$.
Since
$$\frac{1}{2}\d_vF(a)=
   s(u)k_va_v-\sum_{w\in\d_0v}\frac{a_{w^+}}{|b_w|},$$
where
$\d_0v=\d v\cap W_u$,
it follows from vanishing of the derivative and from that
the length function
$a$
is nonnegative that if
$a_v=0$
for some vertex
$v\in u$,
then the function
$a$
vanishes at all vertices of the component
$u$.
Since
$a\not\equiv 0$,
there is a component
$u\in U$,
over which the function
$a$
is positive. Using the BKN-equation for
$(a,\ga)$
we obtain from vanishing of the derivative that
\[\label{eq:zeroangle}
  \sum_{w\in\d v\sm\d_0v}\frac{\ga_w}{|b_w|}a_{w^+}=0 \tag{\rm G}\]
for all
$v\in u$.

Define functions
$a':V\to\R$, $\ga':W\to\{0,\pm 1\}$
by
$a_v'=a_v$
for all
$v\in u$
and
$a_v'=0$
for all
$v\not\in u$;
$\ga_w'=\ga_w$
for all
$w\in W_u$
and
$\ga_w'=0$
for all
$w\not\in W_u$.
It follows from the equality (\ref{eq:zeroangle}) that
$(a',\ga')$
is a compatible symmetric solution to the BKN-equation over
$X$,
and from the preceding we obtain that
$|\ga_w'|=1$
for all edges
$w\in W$
with
$a_{w^-}'\cdot a_{w^+}'\neq 0$.
This is the case (2).

Assume now that
$a_v>0$
for all vertices
$v\in V$.
We put
$a'=a$; $\ga_w'=\ga_w$
for all components
$u\in U$
and all edges
$w\in W_u$
and
$\ga_w'=0$
for all edges
$w$
connecting vertices from different sign components.
Then as above
$(a',\ga')$
is a compatible symmetric solution to the BKN-equation over
$X$
with positive length function,
$a'>0$,
and angle function
$\ga'$
taking the values in
$\{0,\pm 1\}$.
This is the case (3).
\end{proof}

\subsection{Decomposition principle}\label{subsect:decomprin}

Here we prove Lemmas~\ref{Lem:posideal} and
\ref{Lem:idipole}. Our main tool is the decomposition of
labeled graphs
$X=(\Ga,|B|,K)$
into dipoles.

\subsubsection{Dipole}

The graph with two vertices
$v$, $v'$
and one nonoriented edge
$(w,-w)$
between them,
$w\in\d v$,
will be called {\em the dipole.} The BKN-equation over
a labeled dipole
$D=(k_w,k_{-w};b_w)$
is
$$k_wa_w=\frac{\ga_w}{|b_w|}a_{-w}, \quad
   k_{-w}a_{-w}=\frac{\ga_{-w}}{|b_w|}a_w,$$
where
$k_w$, $k_{-w}$
are the charges of the vertices
$v$, $v'$
respectively,
$b_w\in\Z\sm\{0\}$
is the intersection index. We note that labeled dipoles
(as well as other labeled graphs) considered here are not
necessarily associated with graph-manifolds.

All NPC-solutions
$(a_w,a_{-w};\ga_w)$
to this equation can be described as follows. The positive
numbers
$a_w$, $a_{-w}$
are defined up to a common positive factor. It follows
from the condition
$a_w\cdot a_{-w}\neq 0$
that
$\ga_w^2=k_wk_{-w}b_w^2$.
This means, in particular, that for different signs of
the charges,
$k_wk_{-w}<0$,
or for
$k_wk_{-w}b_w^2\ge 1$
the BKN-equation of a dipole has no NPC-solution.
Thus we assume that
$0\le k_wk_{-w}b_w^2<1$.

For
$k_wk_{-w}>0$
the sign of
$\ga_w$
coincides with the sign of the charges
$k_w$, $k_{-w}$
and
\[\label{eq:dipole}
 \left(\frac{a_w}{a_{-w}}\right)^2=\frac{k_{-w}}{k_w}.\tag{\rm dip}\]
In the case
$k_w=k_{-w}=0$,
the equation put no restriction to
$(a_w,a_{-w})$.
When one of the charges vanishes and another not, there obviously is
no NPC-solution.

\subsubsection{Conjunction and decomposition of labelled graphs}

Under conjunction operation of labeled graphs
$X'$, $X''$
(the case
$X'=X''$
is not excluded) some of their vertices are glued and their charges
are added (notation:
$X=X'+X''$).
This operation corresponds to toral summation of oriented
graph-manifolds along some maximal blocks (for details see
\cite{BK1}). The converse operation is called a decomposition
of a labeled graph.

The decomposition principle claims roughly speaking that
compatible solutions to the BKN-equation are parameterized
(with some reservation) by decompositions of the labeled graph
into labeled dipoles satisfying appropriate conditions.

Here is the precise statement. The decomposition principle says that

(a) if
$(a',\ga')$, $(a'',\ga'')$
are compatible solutions to the BKN-equations over
labeled graphs
$X'$, $X''$
respectively,
$X=X'+X''$
is the conjunction along some vertices, and the values of
the functions
$a'$, $a''$
at the glued vertices coincide (in the case
$X'=X''$
we require
$(a',\ga')=(a'',\ga'')$),
then
$(a,\ga)$
is a compatible solution to the BKN-equation over the labeled graph
$X$,
where the functions
$a$, $\ga$
are naturally defined by the solutions
$(a',\ga')$, $(a'',\ga'')$:
their values coincide with the values of corresponding functions
of that solution which is defined along the considered
vertices or edges.

(b) Conversely, any compatible solution
$(a,\ga)$
to the BKN-equation over a labeled graph
$X$
with positive length function
$a$
canonically defines a decomposition of the graph
$X=\sum_{e\in E}D_e$
into labeled dipoles
$D_e=(k_w,k_{-w};b_w)$,
where
$e=(w,-w)$
is a nonoriented edge, according to the conditions:
their intersection indices
$b_w$
are the same as of the graph
$X$;
for every edge
$w\in W$
the charge of the vertex
$v=w^-$
of the decomposition dipole
$D_{(w,-w)}$
is defined by
$k_w=\frac{\ga_w}{|b_w|}\cdot\frac{a_{w^+}}{a_{w^-}}$.
Furthermore, the solution
$(a,\ga)$
restricted to any decomposition dipole is a compatible
solution to the BKN-equation for that dipole.

Here we use the notation of an edge
$w$
as the subscript in the notation of the vertex charge of a dipole
because this does not lead to any ambiguity and distinguishes the obtained
charge from the corresponding charge of
$X$.
Moreover, this is convenient for writing the relation
$k_v=\sum_{w\in\d v}k_w$, $v\in V$,
between the charge
$k_v$
of
$X$
and the charges
$k_w$
of the dipoles into which
$X$
is decomposed by a compatible solution.

\subsubsection{The vertex balance}
The mentioned above equality
$k_v=\sum_{w\in\d v}k_w$, $v\in V$,
we call the equation of the vertex balance. It is a necessary
condition for the decomposition of a labeled graph into dipoles
defined by a compatible solution to the BKN-equation.
Moreover, the vertex balance equation together with
the BKN-equations over the dipoles of a decomposition
$X=\sum_{e\in E}D_e$
for collections
$(a_w,a_{-w};\ga_w,\ga_{-w})$, $e=(w,-w)$,
implies the BKN-equation over
$X$
for
$(a,\ga)$,
if it is possible to adjust the collections
$(a_w,a_{-w})$
by positive factors so that their elements coincide at
the glued vertices and hence define a required length function
$a$.
In particular, if the graph
$\Ga$
of a labeled graph
$X=(\Ga,|B|,K)$
is a tree, then any decomposition
$X=\sum_{e\in E}D_e$
into labeled dipoles, each of which
$D_e=(k_w,k_{-w};b_w)$, $e=(w,-w)$,
has a compatible symmetric solution
$(a_w,a_{-w};\ga_w)$
to its BKN-equation with positive
$a_w$, $a_{-w}$,
defines a compatible symmetric solution
$(a,\ga)$
to the BKN-equation over
$X$.
Namely,
$\ga(w)=\ga_w$
and
$a_v/a_{v'}=a_w/a_{-w}$
for every edge
$w\in W$,
where
$v=w^-$
is the beginning of
$w$,
and
$v'=w^+$
is its end. Since there is no nontrivial circuit, it is
possible to define inductively a function
$a:V\to\R$
using the relations
$a_v/a_{v'}=a_w/a_{-w}$,
the restriction of which to every dipole
$D_e$
differs from
$(a_w,a_{-w})$, $e=(w,-w)$,
possibly by a common positive factor only. Then
$(a,\ga)$
is a required solution to the BKN-equation according (a)
of the decomposition principle.

\subsubsection{Polarized discharge of a vertex}

\begin{Lem}\label{Lem:poldischarge} Let
$X=(\Ga,|B|,K)$
be a labeled graph with simply connected graph
$\Ga$.
Given its vertex
$v\in V$,
there are
$\de\in\{0,\pm 1\}$
and
$\ep>0$
such that if
$|k_v'|<\ep$, $\sgn(k_v')=\de$,
and the charges of the remaining vertices are the same as of
$X$, $k_{v'}'=k_{v'}$,
then the BKN-equation over the labeled graph
$X'=(\Ga,|B|,K')$
has a NPC-solution.
\end{Lem}

\begin{proof} We argue by induction over the number
of vertices. For a dipole the assertion immediately follows
from the description of NPC-solutions to its BKN-equation.
In general case, if
$v$
is an interior vertex of the tree
$\Ga$,
then
$X$
can be represented as the conjunction along
$v$
of a finite collection of trees each of which has less vertices
than
$X$.
Then the assertion follows from the inductive assumption and (a)
of the decomposition principle because the graph
$\Ga$
is a tree.

Now, assume that for the vertex
$v$
there is a unique neighboring vertex
$v'$
in
$\Ga$.
Then
$X$
can be represented as the conjunction along
$v'$
of a dipole
$D$
with vertices
$v$, $v'$
and a labeled tree
$X_1$.
By the inductive assumption, one can suppose that
the charge
$k_{v'}^1$
of the vertex
$v'$
of
$X_1$
is close to zero and has a sign so that the BKN-equation over
$X_1$
has a NPC-solution. Then the dipole charge
$k_{v'}^D$
is defined from the vertex balance
$k_{v'}^D+k_{v'}^1=k_{v'}$.
However, whichever is the charge
$k_{v'}^D$,
it is always possible to find a sign
$\de\in\{0,\pm 1\}$
so that if the dipole charge
$k_v^D$
at
$v$
has the sign
$\de$
and is sufficiently close to zero, then the BKN-equation over
$D$
has a NPC-solution. As above, the assertion now follows from
the decomposition principle.
\end{proof}

For the proof of Lemmas~\ref{Lem:posideal} and
\ref{Lem:idipole} we shall use the following characteristic of ideals.

\begin{Lem}\label{Lem:ideal} Suppose that a labeled graph
$X=(\Ga,|B|,K)$
has the following property: for every charge collection
$K'\in\R^V$
from some neighborhood in
$\R^V$
of
$K$,
the BKN-equation over
$X'=(\Ga,|B|,K')$
has a NPC-solution. Then
$X$
is an ideal.
\end{Lem}

\begin{proof} Assume that
$X$
is the subgraph of a connected labeled graph
$X_0$.
We describe how to find a NPC-solution to the BKN-equation over
$X_0$.
The graph
$X_0$
can be represented as the conjunction of
$X'$,
a disjoin collection of trees
$Y$
and a collection of dipoles
$Z$
with zero charges,
$$X_0=X'+Y+Z.$$
Furthermore,
$X'=(\Ga,|B|,K')$,
and every tree
$Y_{\al}$
from
$Y$
has a charge of the conjunction vertex which is sufficiently
close to zero and satisfies the condition of
Lemma~\ref{Lem:poldischarge}. It follows from the vertex balance
that the charge collection
$K'$
can be chosen arbitrarily close to the initial collection
$K\in\R^V$.
Since the dipoles with zero charges put no restriction to
the length function of a solution to the BKN-equation,
we complete the proof applying the decomposition principle
and taking into account that the collection
$Y$
consists of trees.
\end{proof}

\begin{proof}[Proof of Lemma~\ref{Lem:posideal}] To be definite,
we assume that all charges of the labeled graph
$X$
are positive. Arguing by induction over the number of vertices,
we assume that the assertion is true for all labeled graph with
number of the vertices less than that of
$X=(\Ga,|B|,K)$.
By the condition, the quadratic form
$A_X^+$
on
$\R^V$
has a negative eigenvalue. There is a sufficiently small
neighborhood in
$\R^V$
of
$K$
such that for every
$K'$
from this neighborhood the form
$A_{X'}^+$, $X'=(\Ga,|B|,K')$,
has a negative eigenvalue.

Consider the family of forms
$A_{X',t}^+$, $0\le t\le 1$,
$$(A_{X',t}^+x,x):=\sum_{v\in V}k_v'x_v^2
  -\sum_{w\in W}\frac{t}{|b_w|}x_{w^-}x_{w^+}.$$
For
$t=0$
the corresponding form of this family is positive definite, and for
$t=1$
it has a negative eigenvalue. Thus there is
$t_0\in(0,1)$
such that the form
$A_{X',t_0}^+$
has a zero eigenvalue. Let
$a_0\in\R^V$
be a corresponding eigenvector. Changing the sign of
$\ga_w^0=t_0$
for the edges
$w$
leaving the vertices at which
$a_0$
is negative and simultaneously changing the sign of the negative
values of
$a_0$,
one can assume that
$a_0\ge 0$.
Then
$(a_0,\ga_0)$
is a symmetric solution to the BKN-equation over
$X'$,
for which the length function
$a_0$
is nonnegative, and the angle function
$|\ga_0|<1$.
One can also assume that this solution is compatible.

Consider the labeled subgraph
$X_0\sub X'$
spanned by the vertices at which
$a_0$
is positive. The restriction of
$(a_0,\ga_0)$
to this subgraph is a NPC-solution to the BKN-equation over
$X_0$
(the subgraph has, obviously, no isolated vertex, i.e., a vertex
without incident edges). It follows that
$(A_{X_0}^+a_0,a_0)<0$
for the form
$A_{X_0}^+$.
Thus this form has a negative eigenvalue.

If
$X_0$
is a proper subgraph in
$X'$,
then by the inductive assumption it is an ideal. Thus in any
case,
$X_0=X'$
or
$X_0\neq X'$,
the BKN-equation over
$X'$
has a NPC-solution. By Lemma~\ref{Lem:ideal}, the labeled graph
$X$
is an ideal.
\end{proof}

\begin{proof}[Proof of Lemma~\ref{Lem:idipole}]
In the case of a linear
$\Ga$,
every labeled graph
$X'=(\Ga,|B|,K')$
can be represented as the conjunction of
$p$
labeled dipoles
$D_e'=(k_w',k_{-w}';|b_w|)$,
where
$e=(w,-w)$
is a nonoriented edge, satisfying the NPC-condition
$0<\ga_w^2:=k_w'k_{-w}'b_w^2<1$
(see sect.~\ref{subsect:decomprin}), if the charge collection
$K'\in\R^V$
is sufficiently close to
$K$.
By the decomposition principle, the BKN-equation over
$X'$
has a NPC-solution because the graph
$\Ga$
is simply connected. Now, it follows from Lemma~\ref{Lem:ideal}
that
$X$
is an ideal.

In the loop case
($p=1$),
the BKN-equation over
$X$
has a NPC-solution due to the condition
$|k\cdot b|<2$.
This condition is open in the charge space, thus
$X$
is an ideal.

Now, we assume that
$\Ga$
is a circle with odd number
$p\ge 3$
of edges. If every charge of
$X$
is zero, then the functions
$a\equiv 1$, $\ga\equiv 0$
make up a NPC-solution
$(a,\ga)$
to the BKN-equation. Let
$X'$
be a perturbed graph. One can assume that there is a vertex
$v_1\in V$,
with nonzero charge
$k_1'$.
We assume further that the value
$|k_1'|$
is maximal over all vertices. Decompose
$X'$
into dipoles
$D_e'=(k_w',k_{-w}';|b_w|)$
with all nonzero charges so that together with
NPC-condition above the following conditions are fulfilled:
$k_w'=k_{-w}'$
for all edges
$e=(w,-w)$
not incident to
$v_1$.
For the edges
$e_1=(w_1,-w_1)$, $e_p=(w_p,-w_p)$,
incident to
$v_1$, $v_1=w_1^-=w_p^+$,
we require that
$$\frac{k_{w_1}'}{k_{-w_1}'}=\frac{k_{-w_p}'}{k_{w_p}'}.$$
In view of the relation (\ref{eq:dipole}) from
sect.~\ref{subsect:decomprin}, these conditions provide the
existence of a NPC-solution to the BKN-equation over
$X'$.
It is easy to see that if
$K'$
is sufficiently close to
$K$,
then a required decomposition of the labeled circle
$X'$
into dipoles does exist because
$p\ge 3$
is odd. Hence,
$X'$
is an ideal.
\end{proof}

\subsubsection{Historical remarks} The notion of ideals,
the decomposition principle and the lemma about a polarized discharge
of a vertex have appeared in \cite{BK1}, Lemma~\ref{Lem:ideal}
is taken from \cite{BK2}.

\subsection{Symmetric rational solutions}\label{subsect:symratap}

Here, under some restrictions, we find a rational approximation
of NPC-solutions to the BKN-equation.

\begin{Pro}\label{Pro:symratappr} Let
$X=(\Ga,|B|,K)$
be the labeled graph of a manifold
$M\in\fM$
such that if the graph
$\Ga$
is bipartite, then there are at least four vertices
lying in one and the same part with equal nonzero charges.
Then every NPC-solution
$(a,\ga)$
to the BKN-equation over
$X$
can be approximated by rational NPC-solutions.
\end{Pro}

\begin{Rem}\label{Rem:irrdipole} There are graph-manifolds
$M\in\fM$,
for which the BKN-equation has a real NPC-solution
$(a,\ga)$,
and has no rational NPC-solution. As such an
$M$
one can take an``irrational dipole'' whose graph consists of
two vertices
$v$, $v'$
and one (nonoriented) edge
$e$
between them, the charges and the intersection index
satisfy the condition
$0<k_vk_{v'}b_e^2<1$
and
$\sqrt{k_vk_{v'}}\not\in\Q$
(at least one of the blocks of
$M$
with these properties must have singular fibers).

On the other hand, if a manifold
$M$
has a nonzero charge, then it is easy to find its
characteristic covering
$\wt M\to M$
such that the (connected)
$\wt M$
satisfies the conditions of Proposition~\ref{Pro:symratappr}.
In the case when all the charges are equals zero, the functions
$a\equiv 1$
and
$\ga\equiv 0$
form a NPC-solution
$(a,\ga)$
to the BKN-equation as it was already
mentioned in the proof of Lemma~\ref{Lem:negeig}.
Therefore, while considering the virtual properties
($\ve$), ($\vf$),
one can assume that NPC-solutions are rational by
taking an appropriate characteristic covering if necessarily.
\end{Rem}

The proof of Proposition~\ref{Pro:symratappr} is based on
properties of the incidence matrix of the graph
$\Ga$,
which is defined as the map
$I=I(\Ga):V\times E\to\{0,1,2\}$
($E$
is the set of nonoriented edges of
$\Ga$),
$I_{ve}\neq 0$
if and only if the vertex
$v\in V$
and the edge
$e\in E$
are incident, furthermore,
$I_{ve}=2$
for the loops
$e$
with the vertex
$v$
and only for them. The incidence matrix can also be
identified with a linear map
$I:\R^E\to\R^V$
by considering every restriction
$I_v=I|v\times E$
as the row corresponding to the vertex
$v\in V$.

\begin{Lem}\label{Lem:incirang} If a connected graph
$\Ga$
is not bipartite, then the rank of the incidence matrix
$I$
is equal to the number of the vertices,
$|V|$.
If the graph
$\Ga$
is bipartite, then the rank of
$I$
is equal to
$|V|-1$.
\end{Lem}

\begin{proof} Assume that there is a zero linear
combination of rows,
$\sum_{v\in V}c_vI_v=0$.
Then for every edge
$e\in E$,
we have
$c_v+c_{v'}=0$,
where
$v$, $v'$
are the ends of
$e$
(possibly coinciding), because the column of
$I$
corresponding to
$e$
has nonzero elements for the vertices
$v$
and
$v'$
only. It follows that if the graph
$\Ga$
has a circuit with odd number of edges, i.e., it is not
bipartite, then
$c_v=0$
for the vertices of that circuit, and hence for all
vertices
$v\in V$
by connectedness of
$\Ga$.
In this case, the rank of
$I$
is equal to
$|V|$.

Assume now that the graph
$\Ga$
is bipartite, and let
$V=V_0\cup V_1$
be a decomposition of the vertex set into (nonempty)
parts. Then
$$\sum_{v\in V_0}I_v-\sum_{v\in V_1}I_v=0$$
is the unique (up to a nonzero factor) nontrivial
combination of rows of
$I$,
and hence its rank is equal to
$|V|-1$.
\end{proof}

There are several interpretations of the BKN-equation.
The most convenient one for symmetric rational
approximations is the following. For a positive length function
$a$,
the BKN-equation over
$X$
can be written in more symmetric equivalent form
 $$k_va_v^2=\sum_{w\in\d v}\frac{\ga_w}{|b_w|}a_{w^+}a_{w^-},\quad v\in V.$$
If in addition the angle function
$\ga$
is symmetric, then this form in turn is equivalent to the equation
\[\label{eq:symbkn}
   A=I(D),\tag{+}\]
where
$I:\R^E\to\R^V$
is the incidence matrix,
$A\in\R^V$, $A_v=k_va_v^2$,
and
$D\in\R^E$, $D_e=\frac{\ga_w}{|b_w|}a_{w^+}a_{w^-}$
for every edge
$e=(w,-w)$.

\begin{proof}[Proof of Proposition~\ref{Pro:symratappr}: the graph
$\Ga$
is not bipartite]
In this case, we have
$|V|\le|E|$,
because the graph
$\Ga$
is not simply connected. Since the rank of the matrix
$I$
equals
$|V|$,
there is a linearly independent collection of its rows
$E_0\sub E$, $|E_0|=|V|$,
and by a classical theorem from linear algebra,
the solutions to the symmetrized BKN-equation (\ref{eq:symbkn})
can be parameterized by
$|V|$
linear functions
$L_e$, $e\in E_0$
with rational coefficients depending on
$|E|=|V|+(|E|-|V|)$
arguments
$A_v$, $v\in V$, $D_{e'}$, $e'\in E\sm E_0$,
$$\frac{\ga_w}{|b_w|}a_{w^+}a_{w^-}=
  L_e(A_v,D_{e'}),\quad e=(w,-w)\in E_0.$$
Approximating the numbers
$a_v$, $v\in V$,
and
$\ga_{w'}$, $(w',-w')\in E\sm E_0$
by rational numbers
$a_v'$
and
$\ga_{w'}'$
respectively, we obtain a rational approximation
$$\ga_w'=\frac{|b_w|}{a_{w^+}'a_{w^-}'}L_e(A_v',D_{e'}')$$
of
$\ga_w$
for all
$(w,-w)\in E_0$.
This gives a required approximation of the solution
$(a,\ga)$
to the BKN-equation by compatible, symmetric, rational
solutions.
\end{proof}

Now, we assume that the graph
$\Ga$
is bipartite, and let
$V=V_0\cup V_1$
be a decomposition of its vertex sets into parts,
$U\sub V_0$
be a subset consisting of four vertices with equal charges,
$k_v=k\neq 0$
for all
$v\in U$.
Consider the set
$Q=\set{x\in\R^V}{$\sum_{v\in V_0}k_vx_v^2=\sum_{v\in V_1}k_vx_v^2$}$.
The equality
$\sum_{v\in V_0}I_v=\sum_{v\in V_1}I_v$
implies that for any solution
$(a,\ga)$
to the BKN-equation (\ref{eq:symbkn}), the length
function
$a$
lies in
$Q$,
in particular, the set
$Q$
contains points with all coordinates different from zero.

We use the condition
$|U|=4$
to apply a classical result of Lagrange that every natural
number can be represented as the sum of four squares of
integer numbers for the proof of the following Lemma.

\begin{Lem}\label{Lem:ratdensity} The points with rational
coordinates are dense in
$Q$.
\end{Lem}

\begin{proof} W.L.G. we assume that
$k_v\neq 0$
for all
$v\in V$.
Put for brevity
$\ka_v:=-k_v$
for
$v\in V_0$
and
$\ka_v:=k_v$
for
$v\in V_1$.
Then
$Q=\set{x\in\R^V}{$\sum_{v\in V}\ka_vx_v^2=0$}$.

The Lagrange theorem implies that the set
$Q$
contains a point
$c\neq 0$
with rational coordinates. Indeed, approximating a point
$a\in Q$
with
$a_v\neq 0$
for all
$v\in V$
by rational points from
$\R^V$,
we find a rational
$q>0$,
which can be represented as
$$q=\frac{1}{k}\sum_{v\in V\sm U}\ka_vc_v^2,$$
where
$c_v\in\Q$
for all
$v\in V\sm U$.
Representing
$q$
as the sum of four squares of rationals,
$q=\sum_{v\in U}c_v^2$,
by the Lagrange theorem, we get a required point
$c\in Q$.

It suffices to show that the rational points are dense in
the subset
$Q_0\sub Q$
of the points which have all coordinates different from zero.
Take
$b\in Q_0$
assuming that
$c\neq b$.
Moreover, changing if necessarily the signs of coordinates of
$c$,
one can assume that
$A:=\sum_{v\in V}\ka_vb_vc_v>0$.

Consider the line
$\si$
in
$\R^V$
through the points
$c$
and
$b$, $\si(t)=c+t\xi$,
where
$\xi=b-c$, $t\in\R$.
The set
$Q$
is the zero set of the function
$F:\R^V\to\R$, $F(x)=\sum_{v\in V}\ka_vx_v^2$.
The gradient of
$F$
is different from zero at any point
$x\in Q\sm\{0\}$,
thus
$Q\sm\{0\}$
is a regular hypersurface in
$\R^V$.
The restriction of
$F$
to the line
$\si$
is a quadratic polynomial,
$F(t)=F\circ\si(t)$,
which vanishes at
$t=0,1$.
Thus
$F(t)=A't(t-1)$
for some coefficient
$A'\in\R$.
A direct computation shows that
$A'=-2A\neq 0$.
Consequently, the line
$\si$
intersects transversally
$Q$
at
$b$.
Now, approximating the vector
$\xi$
by vectors
$\xi'\in\R^V$
with rational coordinates, we obtain lines
$\si'$, $\si'(t)=c+t\xi'$,
whose intersection points with
$Q$
gives a rational approximation of
$b$.
\end{proof}

\begin{proof}[Proof of Proposition~\ref{Pro:symratappr}:
$\Ga$
is bipartite] The argument for the nonbipartite case
can be applied also here with the following modifications:
a linearly independent collection of rows
$E_0\sub E$
consists of
$|V|-1$
elements, and the length function is chosen from the set
$Q\sub\R^V$, $a\in Q$.
Lemma~\ref{Lem:ratdensity} provides now a required rational
approximation.
\end{proof}

\section{Proof of spectral criteria}\label{sect:proofspcr}

\subsection{Spectral criterion for
($\imm$)=($\hi$)}\label{subsect:proofhi}
Here we prove Theorem~\ref{Thm:specrimmhi}.

\begin{proof}[Only if]
Assume that the BKN-equation over
$M$
has a compatible solution
$(a,\ga)$,
and we can also assume that the length function
is positive,
$a>0$.
By Lemma~\ref{Lem:nonposeigen}, the operator
$A_M^+$
has a nonpositive eigenvalue, and if it has no
negative eigenvalue, then the length function
$a\in\R^V$
lies in the kernel of this operator and
$\ga_w=\sgn(k_{w^-})$
for every edge
$w\in W$.
Since
$\ga_w\cdot\ga_{-w}\neq -1$
by compatibility of the solution, all charges of
$M$
have one and the same sign.
\end{proof}

\begin{proof}[If] First, consider the case when
the operator
$A_M^+$
has no negative eigenvalue, and all charges of
$M$
have one and the same sign. One can assume that all
charges are positive. By the condition, the operator
$A_M^+$
is degenerate. Let
$a\in\R^V$
be a nonzero vector from its kernel. Then the functions
$|a|$
and
$\ga_w=\sgn(a_{w^-})\cdot\sgn(a_{w^+})$, $w\in W$,
are the length function and the angle function respectively of
a compatible solution to the BKN-equation over
$M$.
By Theorem~\ref{Thm:bknlevel}, the manifold
$M$
has the property
($\imm$)=($\hi$).

Now, we assume that the operator
$A_M^+$
has a negative eigenvalue. Then the argument similar to the
proof of Lemma~\ref{Lem:posideal} shows that the deformation
$A_{M,t}^+$, $0\le t\le 1$,
of
$A_M^+$
has for some
$t\in[0,1)$
a zero eigenvalue, which easily give a compatible
solution to the BKN-equation. Thus also in this case
$M$
has the property
($\imm$)=($\hi$).
\end{proof}

\subsection{Spectral criterion for
($\ebd$)}\label{subsect:proofebd}

Here we prove Theorem~\ref{Thm:specrebd}.

\begin{proof}[Only if] Assume a manifold
$M\in\fM$
has the property
($\ebd$).
According Theorem~\ref{Thm:bknlevel}($\ebd$), there is a compatible,
symmetric solution
$(a,\ga)$
to the BKN-equation over
$M$
such that if
$a_{w^-}\cdot a_{w^+}\neq 0$
then
$\ga_w=\pm 1$
for every edge
$w\in W$.
Define a cocycle
$\la:W\to\Z_2$
putting
$\la_w=\ga_w$,
if
$a_{w^-}\cdot a_{w^+}\neq 0$,
and
$\la_w=1$
otherwise. Let us show that the operator
$A_\la$
is weakly degenerate.

Suppose
$a_v\neq 0$
for some
$v\in V$.
Is suffices to show that
$(A_\la a)_v=0$.
We have
$$(A_\al a)_v=k_va_v-\sum_{w\in\d v}\frac{\la_w}{|b_w|}a_{w^+}=
  k_va_v-\sum_{w\in\d v}\frac{\ga_w}{|b_w|}a_{w^+}=0$$
by the BKN-equation, since
$\la_wa_{w^+}=\ga_wa_{w^+}$
for
$w\in\d v$
whatever
$a_{w^+}$
is zero or not.
\end{proof}

\begin{proof}[If] Assume that the operator
$A_\la$
is weakly degenerate for some class
$\la\in H^1(\Ga;\Z_2)$.
Take a nonzero vector
$x\in\R^V$
with property
$(A_\la x)_v=0$
for all vertices
$v\in V$
with
$x_v\neq 0$.
Define a function
$\ga:W\to\{0,\pm 1\}$
by
$\ga_w=\la_w\sgn(x_{w^-}\cdot x_{w^+})$, $w\in W$.
Then
$(|x|,\ga)$
is a compatible, symmetric solution to the BKN-equation
over
$M$
satisfying the condition
($\ebd$)
of Theorem~\ref{Thm:bknlevel}. Thus
$M$
has the property
($\ebd$).
\end{proof}

\subsection{Spectral criterion for
($\ve$)}\label{subsect:proofve}
Here we prove Theorem~\ref{Thm:specrve}.

\begin{proof}[Only if] Assume a manifold
$M\in\fM$
has the property
($\ve$).
According Theorem~\ref{Thm:bknlevel} there is a compatible,
symmetric solution
$(a,\ga)$
to the BKN-equation over
$M$.
Then by Lemma~\ref{Lem:symnonposeigen} the operator
$H_M$
has a nonpositive eigenvalue.
\end{proof}

\begin{proof}[If] Assume the operator
$H_M$
has a nonpositive eigenvalue. According
Theorem~\ref{Thm:bknlevel}($\ve$), it suffices to show that
the BKN-equation over
$M$
has a compatible, symmetric solution. If
$H_M$
has a negative eigenvalue or the function
$s$
from its definition is zero, then this follows from
Lemma~\ref{Lem:negeig}. Thus we assume that the operator
$H_M$
is degenerate and the function
$s$
takes the values in
$\{\pm 1\}$.
Let
$a\in\R^V$
be a nonzero vector from the kernel of
$H_M$.
Then
$$s(u)k_va_v-\sum_{w\in\d v\cap W_u}\frac{a_{w^+}}{|b_w|}=0$$
for all sign components
$u\in U$
and all vertices
$v\in V_u$.
Define
$\ga_w=s(u)\sgn(a_{w^-}\cdot a_{w^+})$
for
$w\in W_u$
and
$\ga_w=0$
for every edge
$w\in W$
connecting vertices from different sign components. Then
$(|a|,\ga)$
is a compatible, symmetric solution to the BKN-equation over
$M$.
\end{proof}

\subsection{Spectral criterion for
($\fib$)}\label{subsect:proofib}

Here we prove Theorem~\ref{Thm:specrfib}.

\begin{proof}[Only if] By Theorem~\ref{Thm:bknlevel}($\fib$),
there is a function
$\ep:V\to\{\pm 1\}$
and a compatible solution
$(a,\ga)$
to the BKN-equation such that
$a_v>0$
for all
$v\in V$
and
$\ga_w=\ep_{w^-}\ep_{w^+}\sgn(b_w)$
for all
$w\in W$.
It suffices to show that the vector
$x=(x_v)\in\R^V$,
defined by
$x_v=\ep_va_v$,
lies in the kernel of the operator
$A_\rho$.
By Lemma~\ref{Lem:cohomdeg} one can assume that
$$J_M^\rho x=
  \sum_{v\in V}\left(\sum_{w\in\d v}\frac{1}{b_w}x_{w^+}\right)v.$$
Thus
$$A_{\rho}x=\sum_{v\in V}\left(k_vx_v-
      \sum_{w\in\d v}\frac{1}{b_w}x_{w^+}\right)v
     =\sum_{v\in V}\ep_v\left(k_va_v-
     \sum_{w\in\d v}\frac{\ga_w}{|b_w|}a_{w^+}\right)v=0.$$
\end{proof}

\begin{proof}[If] Assume the operator
$A_\rho:\R^V\to\R^V$
is supersingular. Let
$a\in\R^V$
be an element from its kernel with all nonzero coordinates,
$a_v\neq 0$, $v\in V$.
Then for the function
$\ep:V\to\{\pm 1\}$, $\ep_v=\sgn(a_v)$,
the functions
$|a|:V\to\R$
and
$\ga:W\to\{\pm 1\}$, $\ga_w=\ep_{w^+}\ep_{w^-}\sgn(b_w)$,
are the length function and the angle function respectively
of a solution to the BKN-equation, which satisfies the condition
($\fib$)
of Theorem~\ref{Thm:bknlevel}. Thus the manifold
$M$
fibers over the circle.
\end{proof}

\subsection{Spectral criterion for
($\vf$)}\label{subsect:proofvf}

Here we prove Theorem~\ref{Thm:specrvf}.

\begin{proof}[Only if] Assume that a manifold
$M\in\fM$
is virtually fibered over the circle. Then by
Theorem~\ref{Thm:bknlevel} there is a compatible,
symmetric solution
$(a,\ga)$
to the BKN-equation over
$M$
with positive length function,
$a>0$.
Assume the operator
$H_M$
has no negative eigenvalue. Then
$(H_Ma,a)\ge 0$.
On the other hand
$(H_Ma,a)\le 0$
by Lemma~\ref{Lem:symnonposeigen}. Thus
$(H_Ma,a)=0$.
Since
$H_M$
is positive semidefinite, the vector
$a$
lies in its kernel and hence
$H_M$
is supersingular.
\end{proof}

\begin{proof}[If] By Theorem~\ref{Thm:bknlevel}($\vf$)
it suffices to show that the BKN-equation over
$M$
has a compatible, symmetric solution with positive
length function. It the operator
$H_M$
has a negative eigenvalue or the function
$s$
from its definition equals zero, then this follows from
Lemma~\ref{Lem:negeig}.

One can assume that the operator
$H_M$
is positive semidefinite and supersingular, and the function
$s$
takes the values in
$\{\pm 1\}$.
Take
$a\in\R^V$
from the kernel of
$H_M$
such that
$a_v\neq 0$
for all
$v\in V$.
Define a function
$\ga:W\to\{0,\pm 1\}$
putting
$\ga_w=0$
for all edges
$w\in W$
connecting vertices from different sign components, and
$\ga_w=s(u)\sgn(a_{w^-}\cdot a_{w^+})$
for all sign components
$u\in U$
and all
$w\in W_u$.
Then
$(|a|,\ga)$
is a compatible, symmetric solution to the BKN-equation over
$M$
with positive length function.
\end{proof}

\subsection{Spectral criterion for
($\npc$)}\label{subsect:proofnpc}

Here we prove Theorem~\ref{Thm:specrnpc}.

\begin{proof}[Only if] Assume
$M\in\fM$
carries a NPC-metric. Then by
Theorem~\ref{Thm:bknlevel}($\npc$) there is a compatible,
symmetric solution
$(a,\ga)$
to the BKN-equation over
$M$
with positive length function,
$a>0$,
and the angle function
$|\ga|<1$.

Assume the operator
$H_M$
has no negative eigenvalue. Thus
$(H_Ma,a)\ge 0$.
On the other hand
$(H_Ma,a)\le 0$
by Lemma~\ref{Lem:symnonposeigen}. Hence
$(H_Ma,a)=0$
and moreover
$s(u)\ga_w=1$
for every sign component
$u\in U$
and every edge
$w\in W_u$.
Together with condition
$|\ga|<1$
this means that every sign component consists of one point,
and
$W_u=\es$.
In this case
$(H_Mx,x)=\sum_{v\in V}s(v)k_vx_v^2$
for every
$x\in\R^V$.

Assume the function
$s$
takes the values in
$\{\pm 1\}$.
Then
$s(v_0)=1$
for some vertex
$v_0\in V$
with positive charge,
$k_{v_0}>0$.
Since
$0=(H_Ma,a)=\sum_{v\in V}s(v)k_va_v^2$,
this means that the operator
$H_M$
has a negative eigenvalue, a contradiction. Thus
$s\equiv 0$.
\end{proof}

\begin{proof}[If] By Theorem~\ref{Thm:bknlevel}($\npc$)
it suffices to prove that the BKN-equation over
$M$
has a NPC-solution. This is proved in
Lemma~\ref{Lem:negeig}.
\end{proof}

\subsubsection{Historical remarks}\label{subsubsect:gap}

As we already mentioned, the spectral criterion for
($\npc$)
obtained in \cite{BK2} is inaccurate. Namely, in terms of this
survey, it claims that the property
($\npc$)
is equivalent to that the operator
$H_M$
has a negative eigenvalue, or
$H_M\equiv 0$.
As a counterexample one can take a graph-manifold
$M\in\fM$,
which is the mapping torus of a Dehn twist of an
orientable surface of the genus
$\ge 2$
along a simple closed curve not homological to zero.
The graph of such a manifold is a loop, the charge of the vertex
equals
$k=2/|b|$,
where the intersection index
$b$
coincides with the order of the twist, and
$H_M\equiv 0$.
However,
$M$
admits no NPC-metric because
$s\not\equiv0$.
Nevertheless, all results of
 \cite{BK2}, which use the mentioned criterion, remain true,
in particular, this holds for the criterion for mapping tori of
Dehn twists to have a NPC-metric (Theorem~4.6.3).

A similar inaccuracy has occurred in the criterion for the
property of an infinite graph-manifold to carry
a NPC-metric with bounded sectional curvatures and finite volume,
obtained in \cite{BK3}. The error is also not a principal one, and it
can be corrected similarly (in term of this survey) by replacing
the condition
$H_M\equiv 0$
by the condition
$s\equiv 0$.


\bigskip

\noindent
St.-Petersburg Dept. of Steklov Math. Institute

\noindent
Fontanka 27,

\noindent
191023, St.-Petersburg, Russia

\noindent
{\tt buyalo@pdmi.ras.ru}

\noindent
{\tt svetlov@pdmi.ras.ru}

\end{document}